\documentclass[11pt]{article}
\usepackage{epsfig,color}
\usepackage{url}
\usepackage[numbers,square]{natbib}

\usepackage{amsfonts, dsfont}
\usepackage{amsmath}
\usepackage{amssymb}
\usepackage{soul}
\usepackage{xcolor}
\usepackage{graphics}
\usepackage{subfigure}
\usepackage{breakurl}

\newcommand{\E}{\mbox{\bf E}}
\newtheorem{Theorem_}{Theorem}
\newtheorem{Proposition_}{Proposition}

\newtheorem{Remark_}{Remark}
\newtheorem{Lemma_}{Lemma}
\newtheorem{Example_}{Example}

\newtheorem{Assumption_}{Assumption}

\newcommand{\commentout}[1]{}
\newcommand{\be}{\begin{equation}}
\newcommand{\ee}{\end{equation}}
\newcommand{\bea}{\begin{eqnarray}}
\newcommand{\eea}{\end{eqnarray}}
\newcommand{\bean}{\begin{eqnarray*}}
\newcommand{\eean}{\end{eqnarray*}}

\def\s{c_p}
\def\h{c_h}
\def\e{c_e}

\newenvironment{pf}{\vspace{0.1in}\noindent \textbf{Proof $\!\!\!$} ~}{~
\hspace*{\fill} $\square$ }

\begin{document}




\begin{center}
{\Large Dynamic Pricing 
in a Dual Market Environment}
\bigskip

\medskip

{Wen (Wendy) Chen}

{\small Providence Business School, Providence College, Providence, RI 02908,\\
wchen@providence.edu}

{Adam Fleischhacker}

{\small Lerner College of Business and Economics, University of Delaware, Newark, DE 19716,\\
ajf@udel.edu}

{Michael N. Katehakis}

{\small Department of Management Science and Information Systems, 
Rutgers University, NJ 08854,\\
mnk@rutgers.edu}

\medskip

\end{center}

\begin{abstract}
This paper is concerned with  the determination of pricing strategies for a  firm that in each period of a finite horizon   receives replenishment  quantities of a single product which it sells in  two   markets,
e.g., a long-distance market  and an on-site market. The key difference between the two markets is that the long-distance market provides for a one period delay in demand fulfillment.
In contrast, on-site orders must be filled immediately as the customer is at the physical on-site location.  We model the  demands in consecutive periods  as   independent random variables and  their distributions depend on the item's price in accordance with two general stochastic demand functions: additive or multiplicative. 
 The firm uses  a single pool of inventory to fulfill demands from both markets.    
We investigate properties of the structure of the dynamic pricing strategy that maximizes the total  expected discounted  profit over the  finite time horizon, under fixed or controlled replenishment conditions.   Further, we  provide  conditions under which  one market may be the preferred outlet to sale over the other.
\end{abstract}


\medskip

\section{Introduction}\label{S_Introduction}

This paper   investigates the problem of a firm that in each period of a finite horizon  $t=1,\ldots,T$,  
  receives replenishment  quantities of a single product
that it sells through two   markets: i) an on-site  market  (e.g., physical stores)  and ii) a long-distance  market (e.g., an online site).
The firm aims to maximize its total expected discounted revenue over a finite sales horizon of $T$ periods by adjusting the selling prices $p_{i,t}$ at the on-site  market  ($i=s$) and the long-distance  market ($i=\ell$) in each period $t$.  Both markets' price dependent demands are satisfied  using inventory that is held at the on-site  market  location.  When inventory is available, the on-site market's demand is satisfied immediately while the long-distance market's demand is satisfied after a one period delay.  In this way, our scenario mimics one where physical channel customers are fulfilled while in the store and the firm's online/long-distance customers are fulfilled after a given delay (i.e. we assume they are willing to wait one period before product is shipped).  From the firm's perspective, online orders become a deterministic component of the next period's demand and hence, prices can be set for the following period while taking advantage of this information.  Inventory in our model is exogenously determined and we consider it to be a predetermined supply schedule as is often the case with fashion items.  For tractability, we assume excess demands from both markets are fully backlogged at a specified per unit per period penalty cost.    For long-distance customers, a reasonable shipping delay of one period is considered palatable.

  Examples of dual market strategies using a common pool of inventory are becoming more   prevalent.  Large retailers like Nordstrom and Macy's are expanding the role of their stores beyond their traditional role and they     are now using a combination of technology and customer service known as omnichannel fulfillment to avoid stockout-driven lost sales from within their stores.
    In addition to being a shopping destination for local customers, these retailers are transforming stores into online order fulfillment centers and are using store inventory to satisfy an additional market  of online shoppers.  Firms employing this type of dual market strategy have found dramatic improvements in inventory turns and reduced markdowns because of  the enlarged customer base through which store inventory can be sold \citep{retailnet2012,wsj2013}.
  In addition, there is currently a trend  for  U.S. retailers to use their online sales channel and U.S. based inventory to reach customers overseas, e.g.,  in China \citep{wsj2014}.  Partnerships with Alipay, a payment processor closely linked to Alibaba Group Holdings, enable U.S. retailers to overcome both  economic and regulatory hurdles to enable this type of transactions. 
 In addition, Chinese consumers often prefer the reliability and brand authenticity offered by buying directly from U.S. stores, cf. 
\citep{bl}.    
    The reverse shipping direction is also becoming more common with an inventory pool in China being used to reach Chinese customers living abroad or overseas \citep{intRetail2014}.

As serving dual markets has become easier for firms to achieve, smaller firms are also able to pursue dual market strategies.  One such firm, J\&R Music, expanded into online sales after the events of September 11, 2001 led to decreased foot traffic at their New York City store \citep{Butler2005}.  The ability to reach new customers allowed them to offset the reduction in local customers.  Gary's Wine and Marketplace, a local New Jersey wine retailer, attributes their winning of Beverage Dynamic's 2012 Retailer of the Year to their new online presence \citep{Khermouch2012}.  The wine seller now has ``10\% of the business'' coming from its online store that is operated out of the back of its flagship Wayne, NJ retail outlet.

In terms of pricing in dual markets, pricing policies within firms can be either constant across markets or vary by market \citep{liu2006,Huang2009}.  In fashion retail, using pricing differences to intensify demand in online and/or physical markets is commonplace. %

To provide decision support and insight for this trend towards a single pool of inventory being used to satisfy demand from two markets (channels), we investigate the dynamic pricing policies of a retailer using store-level inventory to serve both its physical in-store customers and an additional market of long distance customers to whom product is shipped.  This retailer will use their pricing policies in each market to intensify or reduce expected demand to better match its current and anticipated inventory positions.  Since each market's demand pulls inventory from a common pool, a key distinction between serving these two markets is how inventory can be deployed to fulfill demand without penalty.  For customers using the physical channel, demand is preferably satisfied instantaneously while customers in the online channel are   more willing to wait for product to be shipped and delivered. As such, the optimal pricing policies developed in this paper regulate each market's demand intensity to balance the benefits of delivery postponement offered by on-line sales 
against advantages, such as larger margins or smaller demand variability, that are available through the physical channel.

\begin{figure}[h]
    \centering
    \subfigure[Online Only Promo]{
    \includegraphics[height=0.25\textwidth]{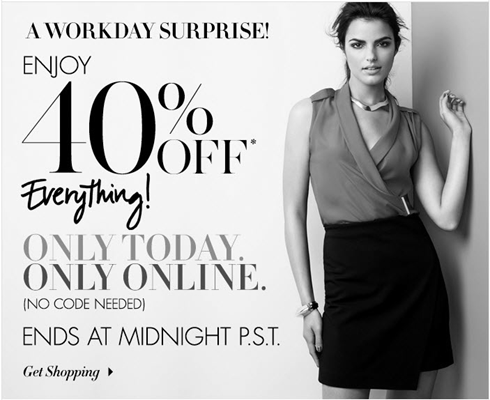}
    \label{fig:subfig1}
    }
    \subfigure[Store Only Promo]{
    \includegraphics[height=0.25\textwidth]{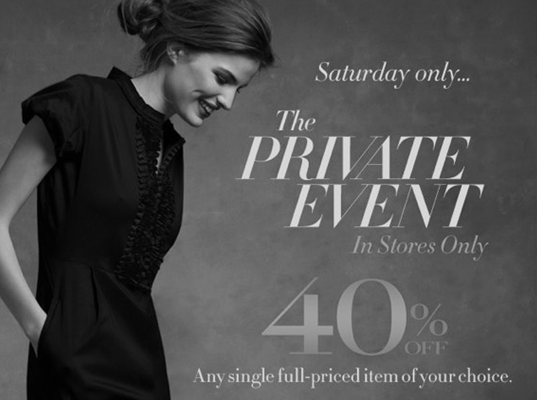}
    \label{fig:subfig2}
    }
    \subfigure[Dual Channel Promo]{
    \includegraphics[height=0.25\textwidth]{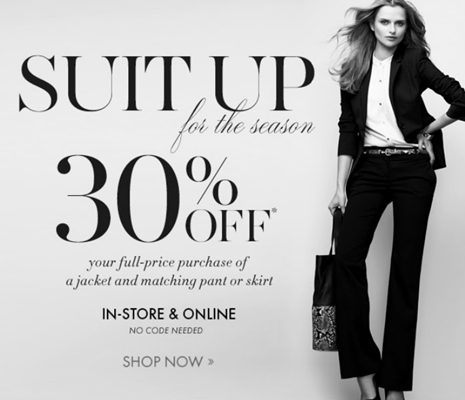}
    \label{fig:subfig3}
    }
   \caption{Retailer Ann Taylor is often seen using promotions to increase demand through its online channel, its physical channel, or both \citep{annTaylor2011}.}
   \label{Fig_AnnTaylor}
\end{figure}

%
%
For example, Ann Taylor, the upscale women's apparel retailer, 
  often employs such market-specific promotional strategies.  In one  promotion it may seek to increase online (long-distance market) demand (e.g. Figure~\ref{fig:subfig1}), with another it seeks to increase in-store (on-site market) demand (e.g. Figure~\ref{fig:subfig2}), and at other times the firm seeks to increase demand in both markets (e.g. Figure~\ref{fig:subfig3}).

%
%
%

The paper has the following structure. In Section \ref{L_Review} we review related
literature. Our study of dynamic pricing extends the existing single market literature to a dual market environment.  %
In Section \ref{S_Model} we   define the notation and state the main assumptions
that we will use throughout the paper.
In Section \ref{S_Discussion} we investigate the basic properties of
model specified in \S\ref{S_Model}. In Section  4, through exploration of both additive and multiplicative cases of demand variability, we find that a different type of demand variability leads to different pricing policies being optimal.
 In \S\ref{SS_Additive}, it is shown  that the optimal prices for both markets  decrease in the inventory level when the
on-site market demand noise is additive. However, this property does not hold in the case at which both markets' demand noise is multiplicative. For the multiplicative case, in \S\ref{SS_Multiplicative}, we show that the firm's market preference under an optimal policy can be specified by a threshold policy.    Special cases for each market's demand being correlated and for the price in each market being equal are also explored in this section.
In Section \ref{S_Discussion2}  insights into the managerial implications of the optimal policies are explored and discussed.   For example, it is
pointed out that at lower inventory levels a firm can prefer sales through the long-distance market even in cases where the marginal profit of the on-site market is significantly higher.
 We   conclude this section with a brief discussion of issues related to
adopting alternative or relaxed assumptions to those made in this paper.

\section{Literature Review}\label{L_Review}

This paper's main contribution is the construction of optimal dynamic pricing policies in a dual market environment and the subsequent development of those policies' managerial implications.  Our study of a second additional market extends the vast literature of dynamic pricing studies that have been done in single market environments \citep[see surveys in][]{bitran2003commissioned,elmaghraby2003dynamic,talluri2005theory} to provide decision support for this relatively new capability of companies to profitably employ  long-distance sales markets; the combination of the internet's global reach with efficiencies in both domestic and international shipping has ushered in the idea of profitably serving local and long-distance markets from the same pool of inventory.

The single market studies most related to our research are finite time horizon models where a firm uses dynamic pricing to intensify or reduce stochastic demand in response to current and future supply availability. Authors Gallego and van Ryzin \citep{GallegoRyzin1994} called this \emph{intensity control} and modeled this using a Poisson arrival process with a price-dependent arrival rate.  A related earlier study considered a model for joint pricing and ordering with an exogenously determined stocking policy \citep{cohen1977joint}, \citep{Thomas1974}.  Our model is a dual market extension to the fully backlogged, periodic review, single market model analyzed by Federgruen and Heching where it is demonstrated how to characterize and compute simultaneous pricing and inventory policies \citep{FedergruenHeching1999}. The work of Chen and Simchi-Levi \citep{ChenSimchi-Levi2004a} extends \citep{FedergruenHeching1999} to model an additional fixed cost component to ordering costs as well as more general demand processes.  Other extensions include assumptions to accommodate substitute products \citep{dong2009dynamic}, incorporation of lost sales as opposed to backlogging \citep{song2009technical}, and the analysis of demand learning with finite capacity \citep{araman2009dynamic}.

Chen in \citep{Chen2001}   provides a model of a firm that may segment its market to gain advanced demand information,  as follows. The customer population is assumed to consist of M segments or types, characterized by different reservation prices Poisson arrival processes with
different rates. Arriving customers are presented with a price schedule that specifies decreasing prices a customer will pay if he/she agrees to different, increasing, shipping delays.  Under sufficient assumptions  
it is  shown  how to  assign 
  customers  to the different market segments so as to maximize the firm's long-run average profits.   In addition  in [21] given an optimal price schedule,  
the author develops a optimal replenishment policy, for a   
a firm that operates with  an $N$-stage supply chain; where  Stage 1 (the     point from which the product is shipped to customers) 
 is replenished by Stage 2, which is replenished by  Stage 3, etc., and Stage N by an outside supplier with ample stock. 
  It is shown that for any  price schedule, the optimal replenishment policy (that minimizes 
  long-run average systemwide holding and backorder costs)
  is to follow an echelon base-stock policy with order-up-to a level that is a function of stage and time. Heuristic computational methods are also provided.

  Our study herein focuses on a  different problem for a firm in which there is a  finite horizon selling season,  and  prices are 
  dynamically adjusted every period as a function of the  ``current'' inventory level and the number of periods remaining in the season. Other  significant differences between the two studies include  the demand models we use and the exogenously determined replenishments in our model.

In extending the dynamic pricing literature to dual markets, we adopt   assumptions that have been used in previous single market dynamic pricing studies.  These assumptions include   the modeling of demand, leadtimes, stockouts, and review policies.  Our first notable assumption is that demand is a non-stationary linear function of price which includes either an additive or multiplicative stochastic term; the  same form that  is used in \cite{AgrawalSeshadri2000,ChenSimchi-Levi2004a}.  For the additive type of demand uncertainty,\citep{ChenSimchi-Levi2004a} prove the optimality of (s,S) policies, but also show that those results do not hold for more general demand functions (i.e. multiplicative plus additive). In most similar studies, optimal prices are shown to be decreasing in inventory level \citep[see for example][]{ChenSimchi-Levi2004b,ChenSimchi-Levi2006,Kazaz2004,Petruzzi1999}.  

 Given the complexity of our model, we adopt the simplifying assumption of a one-period delay (lead time) for the fulfillment of the long-distance market demand.  The use of a one period lead time has a long history of being used to facilitate tractability  \citep[e.g.][]{karlin1958inventory,barankin1961delivery,fukuda1964optimal,veinott1966status,wright1968optimal,dave1982probabilistic}.  More recently, the assumption has been used in \citep{arslan2007single} to develop a mechanism of cost evaluation and optimization for a deterministic replenishment lead-time model (also see this paper for   a detailed review of continuous-review inventory models), in \citep{xu2005multi} to analyze the effect of cancellation contracts on buyer ordering and supplier production policies, and in \citep{wang2009inventory} to study a multi-period inventory model in which a supplier provides two alternative lead-time choices to customers, either a short or a long lead time.

Lastly, in regards to our modeling assumptions, we contribute to a long history of studying fully backlogged periodic review inventory systems to yield tractable insights \citep[see][for examples]{AvivFedergruen2001a,AvivFedergruen2001b,FedergruenHeching1999,BernsteinFedergruen2007,BernsteinFedergruen2004}. Other relevant studies employing this assumption include analyzing pricing and inventory decisions when the supply chain includes multiple retail locations \citep{federgruen2002multilocation}, Markovian demand \citep{Yin2007113}, and stochastic leadtimes \citep{pang2012technical}.

Many notable models in the literature leverage different assumptions regarding the interplay of demand, price, and inventory and we include some similar works here.  These include assumptions of inventory-dependent demand \citep{SmithAchabal1998,datta2001} as well as the modeling of perishable inventory where examples include \citep{bitran1998coordinating} (in the context of a retail chain with coordination among its stores), \citep{feng2000continuous,ZhaoZheng2000,MonahanPetruzziZhao2004,levin2008risk} (whose authors develop a model that incorporates a simple risk measure that can be used to control the probability that revenues are below a minimum acceptable level).  Dynamic pricing policies specific to applications in airline seat pricing are also an important application area and readers are encouraged to see \citep{You1999} and \citep{birbil2010tractable} as examples in this space.

Our research is also related to the dual-market research that is widely studied in the information systems (IS) and marketing literature. This work has mainly  focused on the benefits of online sales and how the additional sales channel leads to higher valuations for a firm \citep[see e.g.][]{BrynjolfssonSmith2000,BrynjolfssonHuSmith2003,CampbellFrei2010}. The key drivers for the higher performance are noted to be better information access and lower setup costs.  Other notable work considers the inter-market competition and demand dependencies that may exist when serving two markets \citep[see][]{FormanGhoseGoldfarb2009,BrynjolfssonHuRahman2009,PentinaPeltonHasty2009}.

The contribution of this paper is important as it provides managerial insight while tackling a level of complexity that has been considered difficult in previous efforts.  Specifically, complexity due to the interplay between a retailer's online and physical channels leads to difficulty crafting optimal dual channel strategies,  see  review by \citep{Agatz2008}. In addition, finding the optimal strategy when more than one type of multiplicative uncertainty is considered creates additional complexity \citep[e.g.][]{AnupindiAkella1993,Chen2013} because a firm cannot always optimally increase price in response to decreasing inventory.  Despite this complexity, we can provide structural insights into the optimal pricing policies in two markets and also, provide managerial insights as to the drivers of preferring demand in one market over another.   There are several papers, tangential to our own, that have investigated  different aspects of this complexity.  Similar to our model, \citep{Chen21092008} model endogenously determined demand in a dual channel environment, but as opposed to manipulating price and/or inventory decisions to intensify demand as we do, the authors restrict attention to policies choosing optimal service levels (as measured by product availability in the physical channel and delivery lead time through the online channel).  Similar studies of endogenously determined response times fall under the category of time-based competition  \citep[see for example][]{Li01021992,Yang2007}. For other related work   we refer to
 \citep{Heching2002}, \citep{qr2012}, \citep{martin2013},  \citep{pis2014} and references therein.

\section{The Basic Model}\label{S_Model}

Let  $T$ denote the finite number of periods in the selling horizon and
let  $p_{i,t}$
denote the posted selling price at market  $i=s,\ell$ in period $t=1,\ldots ,T.$

We will use the following price mean demand model. First, we assume that there exists a function $d$ that represents the relation between the selling price $p_{i,t}$ and the   mean amount of demand so that
$d_{i,t}= d (p_{i,t})$ for each period $t$  and market $i$. Thus,
 deciding the selling price $p_{i,t}$ is equivalent to deciding  the mean amount of demand
$d_{i,t}= d (p_{i,t})$ and conversely setting the
mean demand  of market $i$ at period $t$   to $d$ is equivalent deciding a price: $p_{i,t}(d)=d_{i,t}^{-1}(d)$.
Further, it is assumed that
 $p_{i,t}\in [\underline{\underline{p}}, \overline{\overline{p}}]$  where
$\underline{\underline{p}}$  and $\overline{\overline{p}}$  are known finite  lower and upper
bounds for the price  of market $i=s,\ell$ respectively.
Thus,    $\underline{\underline{d}} = d(\overline{\overline{p}})$
 and $\overline{\overline{d}}= d(\underline{\underline{p}})$ are respectively upper and
 lower bounds for  $d_{i,t}$.

The special cases:
 $d(p)=b -c p$ (for $b>0$, $c >0$, $\underline{\underline{d}}=0$, $\underline{\underline{d}}=b$) and
    $d(p)=b  p^{-c}$ (for $b>0$, $c>1$,$\underline{\underline{d}}=0$, $\underline{\underline{d}}=b-\delta$, $\delta\in (0,b)$)
 are well known examples of mean demand models
 in the   literature and the reader is referred to both \citep{Petruzzi1999} and \citep{ChenSimchi-Levi2004a} for further discussion.
 Throughout the paper, we often simplify the demand model notation by dropping the explicit functional relationship between demand and price and simply write $d$ for $d(p)$.

 The expected revenue  $R_{i,t}(d_{i,t})=d_{i,t} p_{i,t}(d_{i,t})$  is assumed to be a strictly increasing concave function and twice-differentiable in $d$.  The implication of this   assumption is that the  retail firm's revenue increases as the mean demand increases, but does so with diminishing returns to  increasing demand, a similar assumption can be found \citep{ChenSimchi-Levi2004a}.

As in \citep{ChenSimchi-Levi2004a}, we consider two types of demand stochasticity (i.e demand noise), additive and multiplicative. While the firm can choose mean demand by setting price accordingly, actual demand  $D_{i,t}$ is a stochastic function that can be expressed as:
 \bea\label{eq:demand}
D_{i,t} = \epsilon_{i,t} d_{i,t}+ \omega_{i,t},  \   i=s,\ell . \label{eqn_DemandFunction}
 \eea
Here, $\epsilon_{i,t}$ is a random variable representing multiplicative demand noise where $\E(\epsilon_{i,t})=1$ and $\epsilon_{i,t} \in [\underline{\underline{\epsilon_i}}, \overline{\overline{\epsilon_i}}] \subset [0, \overline{\overline{\epsilon_i}}]$, for all
$i,\ t$.
The random variables $\omega_{i,t}$ represent  additive demand noise and it is assumed that $\E(\omega_{i,t})=0$, for all
$i,\ t$.
  Using this expression for stochastic demand, products with price sensitive customers will often exhibit demand realizations consistent with the multiplicative uncertainty model whereas products with less price sensitive customers are best modeled using additive demand uncertainty.  To gain insights, we develop most of our results for the more tractable case where each market's demand is independent of the other.  Subsequently, in \S\ref{SS_Perfect Correlated Demand} we relax this assumption. 

In each period, there are  inventory holding and backorder costs, their   sum we denote by $H(x)$. For our analysis, we adopt the common assumption of a holding cost structure: $H(x)=c_h(x)^+ + c_p(-x)^+$ where $(x)^+=max\{0, x\}$. Generally, we call $c_h$ as  the unit holding cost and $c_p$ as the unit shortage cost. Unmet demand is fully backlogged. 

We note that the backlogging  assumption for items of the primary market $s$ is an approximation, as in reality   (i.e., reasonable values of the pertinent costs)  most of the time backlogging will   be incurred only for items demanded in the secondary market $\ell$. Generally, retailers will satisfy backlogged demand prior to replenishing shelves. However,  $H()$ does not distinguish   the priority with which backlogged customers are satisfied i.e., failing to satisfy a new customer or a backlogged customer incurs the same per period penalty. 
In addition even though, replenishment quantity decisions are made exogenously, the holding cost is  relevant to the retailer as it represents a measure of cost due to some items being  
damaged through mishandling, losses due to  theft or other record keeping problems,  as well as the standard opportunity costs.

The sequence of events  at each period $t=1,\ldots,T$,    can be written as follows.
 The inventory level, $I_t$, is reviewed.
ii) A quantity of product  $q_t$ arrives and is made available in the current period, where  $q_t$ is exogenously determined. 
iii)
 The retailer sets the demand intensity level (through pricing changes) by choosing
  the mean demands for the on-site market $d_{s,t}$ and the long-distance market $d_{\ell,t}$.
iv)
 Demand for the on-site market, $D_{s,t}$, is realized, and satisfied immediately up to the extent there is available supply.  Unmet demand is fully backlogged.
v)
Demand for the long distance market, $D_{\ell,t}$ is realized.  The long-distance market  demand, while known, is neither satisfied nor backlogged in this step.
vi)
Holding costs and backorder costs that are a function $H(\cdot)$  of the
on site market ending inventory $I_t+q_t- D_{s,t} $  are incurred:
$$ H(I_t+q_t- D_{s,t}) =c_h (I_t+q_t-D_{s,t} )^{+} + c_p (I_t+q_t-D_{s,t} )^{-} .$$
Note that inventory held to meet the long-distance demand $D_{\ell,t}$   incurs a holding cost charge in this step, since it is not subtracted from the ending inventory in the term $c_h (I_t+q_t-D_{s,t} )^{+}$; it is the price the firm pays for delaying the fulfillment of $D_{\ell,t}$ (via shipments) by one period.
%
vii)
Demand for the long distance market, $D_{\ell,t}$ is processed   and satisfied to the extent inventory is available.  Unmet demand is fully backlogged and not distinguished from any on-site demand that has been backlogged.

Let $R_{i,t}(d_i)$ ($i=s,\ell$) denote the expected revenue  at market  $i=s,\ell$, in period $t=1,\ldots,T$. Since unmet demand is backlogged, we have:
$$R_{i,t}(d_i)=d_i\,p_{it}(d_i)=\E[(\epsilon_{i,t} d_i +\omega_{i,t})p_{it}(d_s)], \ \mbox{$i=s,\ell$}.$$

 The inventory level for period $t+1$ is calculated as follows:
\bea
I_{t+1}=I_t+q_t-D_{s,t}-D_{\ell,t}.
\eea

%

Let $V_t(I)$ be the optimal expected profit from period $t$ to the end of horizon. The terminal condition is
$$V_{T+1}(I)= c_e min\{0, I\} \mbox{ for all $I$}.$$
  The justification of the terminal condition is as follows. At the end of the horizon, i.e., when  $t=T+1$, if there is a shortage, i.e., when $I<0$, then shortage incurs  is assumed to be supplied by an external supplier at the firm's expense, of  $c_e $ per unit of shortage.
    If there is remaining inventory, $I\geq 0$, then any remaining inventory is assumed to have zero salvage value and $V_{T+1}=0$.

Letting $\alpha \in (0,1]$ be a discount factor, the dynamic programming equations can be written as follows:
\bea
V_t(I) &=& \max_{d_s, d_{\ell} \in [\underline{d}, \overline{d}]} J_t (I, d_s, d_{\ell}), \label{eqn_DP_Vt}
\eea
where
\bea
J_t (I, d_s, d_{\ell}) &=& R_{s,t}(d_s) + R_{\ell, t}(d_\ell) - \E H(I + q_t - \epsilon_{s,t} d_s - \omega_{s,t}  )\nonumber\\
&&\, \, \, \, +  \alpha \E V_{t+1}(I + q_t - \epsilon_{s,t} d_s - \omega_{s,t}  -  \epsilon_{\ell,t} d_{\ell} - \omega_{\ell,t} ).
\label{eqn_DP_Jt}
\eea
The retailer starts each period with a given inventory position $I=I_t$ and a predetermined shipment, $q_t$. The retailer's objective is to control the  demand intensities: $(d_{s,t}(I), d_{\ell,t}(I)),$ through price changes, to maximize expected profit over the planning horizon.

 Let $(d_{s,t}^*(I), d_{\ell,t}^*(I))$ be the maximizer of Eq.~(\ref{eqn_DP_Vt}). We will use the convenetion that   if more than one maximizers exist, the retailer chooses the solution with the smallest sum of each market's mean demand.

To avoid trivial cases, we make the following assumption:
\begin{Assumption_}\label{ass_rev} For   $i=s,\ell$, the following are true
\begin{enumerate}
  \item
 $R_{i,t}^{\prime}(\underline{\underline{d}})>\max\{\s,\h\}$.
   \item
   $\e >\max\{R_{s,T}^{\prime}(0), R_{\ell,T}^{\prime}(0)\} -\s.$
\end{enumerate}
 \end{Assumption_}
The first statement of Assumption~\ref{ass_rev} ensures the retailer has incentive to both carry inventory and backlog demand as required by guaranteeing that the marginal revenue is bigger than both the unit holding cost and the unit shortage cost. The second statement of Assumption~\ref{ass_rev} avoids trivial cases where a retailer chooses to backlog demand towards the end of planning horizon because the outside supplier's cost is competitive with the firm's internal manufacturing costs. Due to this high cost of outsourcing supply in the last period, we assume the retailer will not sell products through the long-distant market in the last period ($d_{\ell,T} = 0$).

\section{Structure of Optimal Policies}\label{S_Discussion}

In this section, we analyze the model specified in \S\ref{S_Model}.  In \S\ref{SS_Concavity}, we investigate the basic properties of  the value function and show how concavity guarantees the existence of an optimal solution to the dual-market dynamic pricing model. In \S\ref{SS_Additive}, we show that the optimal prices for both markets  decrease in the inventory level when the
on-site market demand noise is additive. However, this property does not hold in the case at which both markets' demand noise is multiplicative. For the multiplicative case, in \S\ref{SS_Multiplicative}, we show that the firm's market preference under an optimal policy can be specified by a threshold policy.  A resulting insight is the
identification of conditions under which a firm prefers to sell products in the long-distance market.  Special cases for each market's demand being correlated and for the price in each market being equal are also explored in this section.

\subsection{Concavity and Supermodularity Properties}\label{SS_Concavity}

In the following two lemmas, we state structural results regarding the value function $V_t(I) $ of Eq. (\ref{eqn_DP_Vt}) and the profit function  $J_t(I,d_s,d_\ell)$
of Eq. (\ref{eqn_DP_Jt}).
 We will use these results in  the subsequent analysis.   Lemma ~\ref{lem_ConcaveRevenuefunction} establishes the concavity of $V_t(I)$ as well as concavity and  supermodularity properties for  $J_t(I, d_s, d_\ell)$  and thus, it ensures the existence of an optimal solution.

\begin{Lemma_}\label{lem_ConcaveRevenuefunction}
The following statements are true.  
\begin{itemize}
\item[(i)] $V_t(I)> -\infty$.
\item[(ii)] $V_t(I)$ is concave in $I$ and
$J_t(I,d_s,d_\ell)$ is, component-wise,  concave in each of the variables $I,$ $d_s$.
$d_\ell$.
\item[(iii)] $J_t(I,d_s,d_\ell)$ is submodular in $(d_s,d_\ell)$ and supermodular in $(I,d_s)$ and in $(I,d_\ell)$.
\end{itemize}
\end{Lemma_}

The existence of a maximizer $(d^*_{s,t}(I), d^*_{\ell,t}(I))$ follows from the concavity of $J_t(I,d_s,d_\ell)$.  Further the concavity of the  profit function $J_t(I,d_s, d_\ell)$, in $I$, indicates that the marginal profit decreases in the inventory level.

Combining the concavity results of Lemma~\ref{lem_ConcaveRevenuefunction} and \citep{Topkis1998}, we have $J_t(I,d_s,d_\ell)$ is supermodular in $(I,d_s)$
and $(I,d_\ell)$.  That is to say, if one market's mean demand $d_{i,t}$ (i.e. price) is considered fixed, then the optimal mean demand $d_{j,t}^{*}(I)$ in the other market ($j\neq i$) increases as the inventory $I$ level increases.

To find the optimal mean demand levels $(d_{s,t}^*(I), d_{\ell,t}^*(I))$, we note that the first order partial derivatives can be written as follows:

\bea
\frac{\partial J_t}{\partial d_s} &=& R_{s,t}^{\prime}(d_s)  +  \E\big[ \epsilon_{s,t} H^{\prime}( I + q_t - \epsilon_{s,t} d_s - \omega_{s,t}  \big]\nonumber\\
&& \,\,\,\,\,\,   - \alpha \E \big[\epsilon_{s,t} V_{t+1}^{\prime}(I + q_t  - \epsilon_{s,t} d_s - \omega_{s,t}  -  \epsilon_{\ell,t} d_{\ell} - \omega_{\ell,t})\big]\label{eqn_FOC_ds}\\
\frac{\partial J_t}{\partial d_\ell} &=& R_{\ell,t}^{\prime}(d_\ell)  - \alpha \E \big[\epsilon_{\ell,t} V_{t+1}^{\prime}(I + q_t  - \epsilon_{s,t} d_s - \omega_{s,t}  -  \epsilon_{\ell,t} d_{\ell} - \omega_{\ell,t})\big]\label{eqn_FOC_df}
\eea
From the concavity of $R_{s,t}$, $R_{\ell,t}$, $-H$ and $V_t$, we have both $\frac{\partial J_t}{\partial d_s}$ and $\frac{\partial J_t}{\partial d_\ell}$ decrease in $d_s$ and $d_\ell$. Indeed, note that  $H$ is differentiable almost every where. $V$ is differentiable because demand uncertainty is continuous. First order partial derivatives, as above, can be used because $d_{s,t}^*(I)$ and $d_{\ell,t}^*(I)$ are both continuous in $I$.  Given the concavity of $J_t$,  the optimal solution, $(d_{s,t}^*, d_{\ell,t}^*)$,  is the solution to the above first order partial derivatives equations. Furthermore, we can obtain the following lemma:
\begin{Lemma_}\label{lem_GeneralCase_Increasing}
\bea
\nexists \, \delta>0 : d_{s,t}^*(I+\delta) < d_{s,t}^*(I) \mbox{ and } d_{\ell,t}^*(I+\delta) < d_{\ell,t}^*(I) \label{eqn_GeneralCase_Increasing}
\eea
\end{Lemma_}

Lemma ~\ref{lem_GeneralCase_Increasing} implies that at least one of $d_{s,t}^*(I)$ and $d_{\ell,t}^*(I)$ increases in $I$. Given more inventory on hand, the firm will decrease the sales price in at least one of the two markets. In the following section, we will establish the stronger result that both markets' selling prices are decreasing in inventory level for the case of additive demand noise.

\subsection{Additive Demand Noise}\label{SS_Additive}

In this section, we consider the case where at least one of the market's demand distributions are characterized by purely additive demand noise such that either $\forall t \in T : \epsilon_{s,t} \equiv 1$ or $\forall t \in T : \epsilon_{\ell,t} \equiv 1$.
For this case, price changes in the market(s) with additive demand noise result in changing the mean of a market's demand distribution, but not its variability.  As discussed in \cite{AgrawalSeshadri2000},  additive demand noise is typical of well-established products where the effect of pricing changes on store traffic is well understood.  Uncertainty in these cases tends to be limited to forecasting error.  In contrast to what we will see in the next section when the demand uncertainty is multiplicative, a firm's market preference for selling in one market over another will be unchanged in an additive demand uncertainty environment. The main results of this section,   Theorem \ref{lem_Additiveon-siteDemand},  show that in the additive demand uncertainty environment, the firm prefers to sell more products through both markets when the inventory level $I_t$ or the incoming inventory $q_t$ increases.

\begin{Theorem_}\label{lem_Additiveon-siteDemand}
Under the assumptions made and if the demand noise for the on-site  market is additive, $i.e.$, $\epsilon_{s,t}\equiv 1$ for all $t$, then
\begin{itemize}
\item[(i)] $d_{s,t}^*(I)$ increases in $I$.
\item[(ii)] $I - d_{s,t}^*(I)$ increases in $I$.
\item[(iii)] $d_{\ell,t}^*(I)$ increases in $I$.
\end{itemize}
\end{Theorem_}

\begin{Remark_}\label{Remark_Additiveon-siteDemand_monotonIn_qt}
Following similar logic, we can establish the effect of the incoming inventory $q_t$. If the demand noise for the on-site  market is additive, $i.e.$, $\epsilon_{s,t}\equiv 1$ for all $t$, both $d_{s,t}^*(I)$ and $d_{\ell,t}^*(I)$ increase in $q_t$.
\end{Remark_}

\begin{Remark_}
\label{Remark_lem_Additiveon-siteDemand}
Following a similar argument as that of Theorem~\ref{lem_Additiveon-siteDemand}, we can also obtain that both $d_{s,t}^*(I)$ and $d_{l,t}^*(I)$ are increasing in $I$ if the demand uncertainty in the long-distance market is additive, $i.e.$, $\epsilon_{l,t}\equiv 1$ for all $t$.
\end{Remark_}

Theorem \ref{lem_Additiveon-siteDemand}\emph{(i) and (iii)} states that inventory level increases are accompanied by decreases in the optimal selling prices for both markets and the firm seeks to simultaneously increase demand in both markets. We show in \S\ref{SS_Multiplicative} that these results do not hold in the case of multiplicative demand variability. Also under additive demand uncertainty, the desired amount of customers increases less than one unit when inventory increases by one unit (see Theorem \ref{lem_Additiveon-siteDemand}\emph{(ii)}). Since the marginal revenue of an additional unit of sales decreases as demand increases, the expected demand must increase less than one unit when the inventory level is increased by one
unit.

\subsection{Multiplicative Demand Noise}\label{SS_Multiplicative}

In the multiplicative demand variability case (i.e., $\omega_{i,t}=0$ for $i=s,\ell$ and for all $t$) demand uncertainty increases with increasing demand and equivalently, it increases with decreasing price.  In this section, we will numerically demonstrate that this increased uncertainty leads to non-monotone results. Because of this, characterizing the relationship between inventory and optimal demand levels is more challenging.  However, it is still possible to gain insight and in this section, we leverage the existence of threshold policies to characterize: 1)   conditions under which a retailer makes product available for sale in each market (\S\ref{SSS_Threshold}), Theorem 2)   conditions under which one market is considered preferable because product is made available for sale in that market, but not the other (\S\ref{SSS_MktPref}), Theorem 3)   conditions under which one market is considered preferable because its expected demand is higher than the other market which also offers the product for sale (\S\ref{SSS_MktPref}),Theorem 4) the effects of demand correlation on a firm's market preference (\S\ref{SS_Perfect Correlated Demand}), and Theorem 5) the effects of restricting the price to be equal in both markets on a firm's market preference (\S\ref{SS_Unique_Price}).

Before we present the analytical treatment of the multiplicative case, we use a
numerical example to illustrate the behavior of the optimal average demand $(d_{s,t}^*(I),d_{\ell,t}^*(I))$ as functions of the inventory level.  The parameters of this example are specified in Example \ref{Ex_IllessIs} which is presented along with Figure~\ref{Fig_IllessIs}.  Three important characteristics of the multiplicative demand noise case are established with this example:
\begin{enumerate}
  \item The relationship between inventory level and optimal mean demand $(d_{s,t}^{*}(I), d_{\ell,t}^*)$ is non-monotonic. For example, $d_{\ell,1}^*(-1.3)=0.88<0.97=d_{\ell,1}^*(-1.4)$.
  \item There appears to be a threshold inventory level for each market at which below that inventory level, sales are not pursued in that market and above that level, sales are pursued.  When a market's optimal demand level is zero, $d_{i,t}^*=0$, we say that market $i$ is closed (i.e. $D_{i,t}(0)=0$).
  \item As inventory levels increase, it is possible that the firm's preference to sell in one market over the other (as indicated by a higher expected demand level) may reverse.
\end{enumerate}

Note that in the multiplicative case, an increase in the mean demand $d_{i,t}$ in a particular market is accompanied by an undesirable increase in demand variance $\sigma^2(D_{i,t} $) which makes the characterization of pricing decisions difficult to address.  Similar difficulties have been noted by \cite{LiZheng2006} who study an analogous environment with multiplicative yield uncertainty.  In the multiplicative demand environment with two markets, increasing mean demand is not always the proper response to higher levels of inventory.  In some cases, decreasing demand in one market allows for inventory to build (with the new deliveries) and can be used in subsequent periods.

\begin{samepage}
\begin{figure}[!ht]
\begin{Example_}\label{Ex_IllessIs}
$T=2$, $H(I)=2 \max\{I,0\} + 5 \max\{-I,0\}$, $R_{s,t}(d)=(10-0.5d)d$, $R_{\ell,t}(d)=(9-0.5d)d$ $(0\leq d\leq 9)$ and $V_{T+1}(I)=10\min\{I,0\}$. $\varepsilon_{s,t}$ follows a truncated normal distribution with $E(\varepsilon_{s,t})=1$, $\sigma(\varepsilon_{s,t})=0.6$ and $\varepsilon_{s,t}\in(0,2)$. $\varepsilon_{\ell,t}$ follows another truncated normal distribution with $E(\varepsilon_{\ell,t})=1$, $\sigma(\varepsilon_{\ell,t})=0.9$ and $\varepsilon_{\ell,t}\in(0,2)$. The discount factor $\alpha=0.8$. $q_{T}=1$ and $q_{T-1}=2$. Figure (\ref{Fig_IllessIs}) below illustrates the non-monotonic relationship between the optimal mean demands $d_{i,t}^{*}(I)$ for each market and the inventory level.
\end{Example_}
   \centering
   \includegraphics[width=80mm]{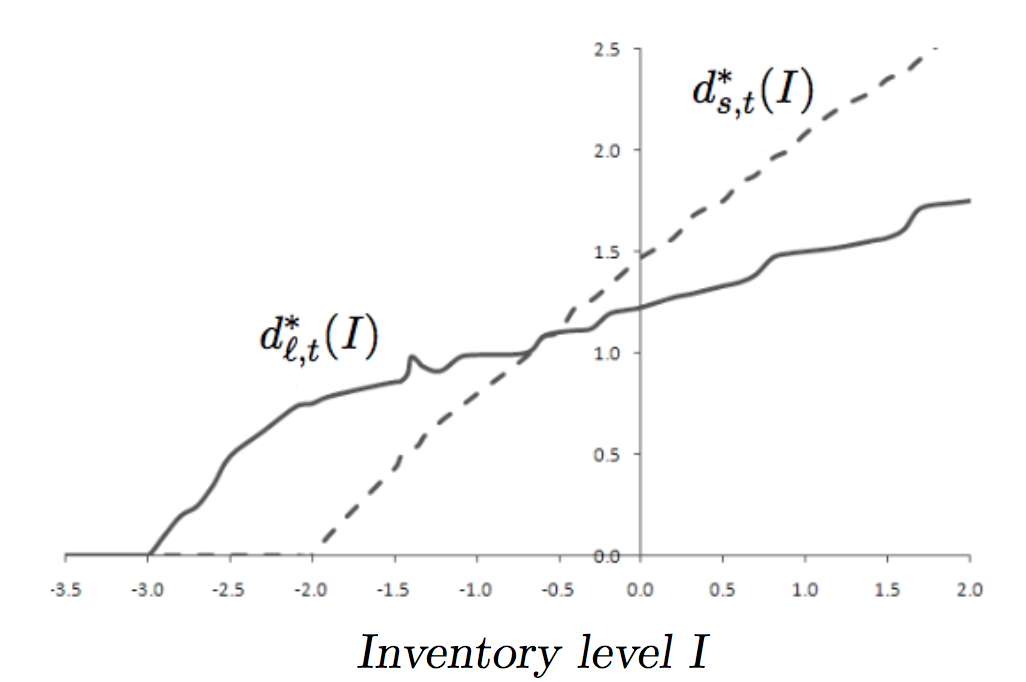}
   \caption{Plot of $d_{s,t}^{*}(I)$ and $d_{\ell,t}^{*}(I)$  versus $I.$}
   \label{Fig_IllessIs}
\end{figure}
\end{samepage}

 \begin{Remark_}\label{Remark_holding}
 
Our insight for the non-monotonicity phenomenon that appears  in the  case with   multiplicative noise, (when
 $D_{i,t} = \epsilon_{i,t} d_{i,t}$ with  $\E(\epsilon_{i,t})=1$) is as follows. When  the controller    increases  the expected demand    $d_{i,t}$ she also increases  its variance with  $d_{i,t}^2$. 
Hence, 
it is possible  that for certain values of inventory increasing $d_{i,t}$ implies 
 cost contributions  of a shortage that are   low relatively to the cost implications  of the additional demand variability. This is phenomenon is further 
amplified by the finite horizon of the selling season. Thus, as   it was also noted in the context of a different model in \cite{FedYang2009},  it may be 
it optimal to target a lower, rather than a higher, expected inventory level after ordering. 
 \end{Remark_}

\subsubsection{Threshold Policies for Market Preferences}\label{SSS_Threshold}
The numerical study shows that optimal selling quantities in on-site market or long distance market may increase or decrease as inventory level increases.  This lack of monotonicity makes analytical insights more difficult to achieve.  
Despite this, we characterize optimal pricing policies using
 the following two simplified benchmark problems:
\bea
\mathbf{(B_\ell)} \, \, \, \, \,  \, \, \, \, \, \, \, \, \, \,  \, \, \, \, \, \, \, \, \, \,  \, \, \, \, \, \, \, \, \, \,  \, \, \, \, \,
V_t^{B_\ell}(I) &=& \max_{d_s, d_\ell} J_t (I, d_s, d_\ell ) \label{eqn_DP_BL_Vt},\\
& & \nonumber\\
J_t^{B_\ell} (I, d_s, d_l ) &=& R_{s,t}(d_s) + R_{\ell,t}(d_\ell) - \E H(I + q_t - \epsilon_{s,t}d_s ) \nonumber\\
    && \, \, \, \, \,  \, \, \, \, \,  + \alpha \E V_{t+1}(I + q_t - \epsilon_{s,t}d_s  - d_\ell)\label{eqn_DP_BL_Jt}.
\eea
\bea
\mathbf{(B_s)} \, \, \, \, \,  \, \, \, \, \, \, \, \, \, \,  \, \, \, \, \, \, \, \, \, \,  \, \, \, \, \, \, \, \, \, \,  \, \, \, \, \,
V_t^{B_s}(I) &=& \max_{d_s, d_\ell} J_t (I, d_s, d_\ell ) \label{eqn_DP_BS_Vt},\\
& & \nonumber\\
J_t^{B_s} (I, d_s, d_l ) &=& R_{s,t}(d_s) + R_{\ell,t}(d_\ell) - \E H(I + q_t - d_s ) \nonumber\\
    && \, \, \, \, \,  \, \, \, \, \,  + \alpha \E V_{t+1}(I + q_t - d_s  - \epsilon_{\ell,t}d_\ell)\label{eqn_DP_BS_Jt}.
\eea
Let $(d_{s,t}^{B_\ell*}(I),d_{\ell,t}^{B_\ell*}(I))$ denote the corresponding maximizer in Eq. (\ref{eqn_DP_BL_Vt}) and $(d_{s,t}^{B_s*}(I),d_{\ell,t}^{B_s*}(I))$ denote the corresponding maximizer in Eq. (\ref{eqn_DP_BS_Vt}) . If there are multiple maximizers, the solution with smallest $(d_{s,t}^{k*}(I),d_{\ell,t}^{k*}(I))$is chosen $(k=B_\ell,B_s)$.

Note that the   above benchmark problems are no longer dynamic programs.  They are one-period representations of our original model (\S\ref{S_Model}) that assume 
 that  $V_{t+1}(\cdot)$ is given and that the demand uncertainty in one of the markets has been eliminated.  In the benchmark problem $B_\ell$, the long-distance market's demand uncertainty in the current period is removed.  Similarly, $B_s$ removes the on-site market's demand uncertainty during the current period.  For these two benchmark problems, threshold policies prescribing the levels of inventory at which the firm makes product available for sale in each market are established below in Lemmas~\ref{lem_BL_BS_Monotone} - \ref{Lem_BL_BS_dlBLanddl}.  These threshold policies are then proven identical to the threshold policies that exist for the original model.  Theorem~\ref{th_MultiplicativeDemand_Threshold}, which is proved with  Lemmas~\ref{lem_BL_BS_Monotone} - \ref{Lem_BL_BS_dlBLanddl}, establishes  this existence and defines threshold inventory levels, $I_{i,t}^*$, such that $d_{i,t}^*(I)=0$ and $d_{i,t}^*(I)>0$ if and only if  $I>I_{i,t}^*$.  We introduce the two benchmark models to yield analytic results by making each market's optimal demand intensity monotonic in the inventory level.  Formally connecting the analytic results of the benchmark problems to the original problem is now shown through Lemmas~\ref{lem_BL_BS_Monotone} - \ref{Lem_BL_BS_dlBLanddl} and Theorem~\ref{th_MultiplicativeDemand_Threshold}.

\begin{Lemma_}\label{lem_BL_BS_Monotone}
\begin{itemize}
\item[(i)] For the benchmark problem ($B_\ell$), we have $d_{\ell,t}^{B_\ell*}(I)$ increases in $I$.
\item[(ii)] For the benchmark problem ($B_s$), we have $d_{s,t}^{B_s*}(I)$ increases in $I$.
\end{itemize}
\end{Lemma_}

\begin{Lemma_}\label{Lem_BL_BS_dlBLanddl}
\begin{itemize}
\item[(i)] $d_{\ell,t}^*(I)=0$ if and only if $d_{\ell,t}^{B_\ell*}(I)=0$.
\item[(ii)] $d_{s,t}^*(I)=0$ if and only if $d_{s,t}^{B_s*}(I)=0$.
\end{itemize}
\end{Lemma_}
%

\begin{Theorem_}\label{th_MultiplicativeDemand_Threshold}
If the demand noise for both  markets is multiplicative, $i.e.$, $\omega_{i,t}=0$ for all $i=s,\ell, \ t=1,\ldots$,T, then the following are true.
\begin{enumerate}
\item[i)]
There exist numbers   $I_{i,t}^*$\, , $i=s,\ \ell$, such that:
\bea
d_{i,t}^*(I) =
\left\{
\begin{array}{cc}
>0 &   \mbox{ if } \,\,\, I>I_{i,t}^*\ ,\\
=0 & \mbox{ if } \,\,\,  I\leq I_{i,t}^* \ .
\end{array}
\right. \label{eqn_MultiplicativeDemand_Threshold}
\eea
\item[ii)] Further, $I_{s,t}^*$ satisfies:
\bea\label{eqn_Marginal_s}
V_t^{\prime}(I_{s,t}^*)= R_{s,t}^{\prime}(0) \, .
 \eea
\end{enumerate}
\end{Theorem_}
Theorem~\ref{th_MultiplicativeDemand_Threshold} proves the existence of these threshold inventory levels, $I_{i,t}^*$, such that $d_{i,t}^*(I)=0$ and $d_{i,t}^*(I)>0$) if and only if  $I>I_{i,t}^*$.  The theorem also provides a  simple characterization  of the on-site market's threshold level $I_{s,t}^*$  in terms of marginal changes of the value function, $V_t(I)$, and the revenue function, $R_{s,t}(d)$.

Mathematically, choosing which products to sell through which combination of markets  is a more challenging problem \citep{zhang2010crafting}.  Retailers sometimes choose to offer greater assortment through their online channel while other retailers have products tailored to local needs that are not available online. The existence of threshold policies shows  that inventory considerations play an important role in this decision. Further characterization of these policies is pursued in subsequent sections.

\subsubsection{Market Preference} \label{SSS_MktPref}
For any time period $t$, we call the market with the higher demand intensity for a given inventory level, $i : \operatorname*{arg\,max}_{i \in \{s,\ell\}} d_{i,t}^*(I)$, as the ``\emph{the preferred market}''.  At lower inventory levels, the preferred market can also be characterized as the market with the lower threshold inventory as the firm will open only the preferred market and keep the mean demand of the other market at $0$.  To differentiate, the market with the lower threshold inventory will be called the ``\emph{preferred opening market}''. As seen in Example \ref{Ex_IllessIs}, market preference may be surprising.  In that example, the firm prefers the long-distance market at lower inventory levels even though the on-site market has greater  expected  revenue ($R_{s,t}(d)>R_{\ell,t}(d)$, for all $d$) and lower demand variability ($\sigma(\varepsilon_{\ell,t})>\sigma(\varepsilon_{s,t})$).   Also from that example, we observe that as the inventory level increases, this preference can change.  For larger values of $I$, market preference is reversed and the on-site  market is preferred as more demand is encouraged through that market, $d^{*}_{s,1}(I) >d^{*}_{\ell,1}(I) $.

We now develop analytic insights into market preferences. From Theorem \ref{th_MultiplicativeDemand_Threshold}, we know that  the
marginal profit of the available inventory equals the marginal revenue of opening the on-site market $R^{'}_{s,t}(0)$ at the on-site market threshold inventory level  $I^*_{s,t}$.   Lemma \ref{lem_MultiplicativeDemand_ChooseingOnline_Threshold_Vprime}
below provides an analogous, albeit complex, expression for the long-distance
market threshhold $I_{\ell,t}^*$.

\begin{Lemma_}\label{lem_MultiplicativeDemand_ChooseingOnline_Threshold_Vprime}
\bea
V_t^{\prime}(I_{\ell,t}^*)=R_{\ell,t}^{\prime}(0) - \E H^{\prime}(I_{\ell,t}^* + q_t - \epsilon_{s,t} d_{s,t}^*(I_{\ell,t}^*))
\label{eqn_MultiplicativeDemand_ChooseingOnline_Threshold_Vprime}
\eea
\end{Lemma_}

The conditions under which one market is preferred to another when the inventory level is low can now be stated:

\begin{Theorem_}\label{lem_MultiplicativeDemand_IsIl} The following two statements are true.
\begin{itemize}
  \item[i)] If $R_{s,t}^{\prime}(0) > R_{\ell,t}^{\prime}(0) +\s$, then $I_{s,t}^*<I_{\ell,t}^*$.
  \item[ii)] If $R_{\ell,t}^{\prime}(0) > R_{s,t}^{\prime}(0) +\h $, then $I_{s,t}^*>I_{\ell,t}^*$.
\end{itemize}
\end{Theorem_}

Managerially, Theorem ~\ref{lem_MultiplicativeDemand_IsIl} provides conditions for having a clearly preferred opening market. In general, if one market's opening has an associated marginal expected revenue that sufficiently exceeds the other market's opening marginal  expected revenue, then that market becomes the preferred market.  The question is what does it mean to sufficiently exceed the other market's marginal revenue?  For the on-site market to be clearly preferred, then the on-site's marginal revenue must exceed the long distance market's revenue by one unit of shortage cost.  For the long distance market to be clearly preferred, then the long distance market's marginal revenue must exceed the on-site market's revenue by one unit of holding cost.  In the next theorem, conditions where the preference is less clear are explored:

\begin{Theorem_}\label{lem_MultiplicativeDemand_preferIl}
When $(R_{\ell,t}^{\prime}(0)-\alpha R_{\ell,t+1}^{\prime}(0))/\alpha>c_p$, the following two statements are true.
\begin{itemize}
  \item[i)] If $R_{s,t}^{\prime}(0) > R_{\ell,t}^{\prime}(0) + c_p$, then $I_{s,t}^*<I_{\ell,t}^*$.
  \item[ii)] If $R_{s,t}^{\prime}(0) < R_{\ell,t}^{\prime}(0) + c_p$, then $I_{s,t}^*>I_{\ell,t}^*$.
\end{itemize}
\end{Theorem_}

Theorem~\ref{lem_MultiplicativeDemand_preferIl} is best interpreted by examining the effects of shortage costs on the  opening market preference. As   in Theorem ~\ref{lem_MultiplicativeDemand_IsIl}, when shortage costs are small and the on-site market provides greater marginal revenue upon market opening, then the on-site market is preferred (Theorem~\ref{lem_MultiplicativeDemand_preferIl}(i)). This preference continues until $c_p$ is large enough to meet the conditions of Theorem~\ref{lem_MultiplicativeDemand_preferIl}(ii).  At this point, the value of delayed fulfillment in the long distance market comes to bear. By fulfilling long-distance orders in the subsequent period, pricing decisions in that subsequent period can be made with less demand uncertainty.  And hence, by selling through the long-distance market the firm has less risk of incurring shortage costs than if the on-site market were opened.  Notice that the firm still prefers selling in the long-distance market when $R_{\ell,t}^{\prime}(0) < R_s^{\prime}(0) < R_{\ell,t}^{\prime}(0) + \s$. That is to say, even though the marginal expected revenue of selling in the long-distance market  is less than the marginal expected revenue of selling in the on-site market, the firm still prefers selling  in the long-distance market.  We again refer the reader to Example~\ref{Ex_IllessIs} where the ``counter-intuitive'' behavior implied by Theorem~\ref{lem_MultiplicativeDemand_preferIl}\emph{(ii)} is seen.  In this example, the on-site market has higher marginal profit (i.e. $R_{\ell,t}(d)<R_{s,t}(d)$) and also lower variability in demand (i.e. $\sigma(\varepsilon_{\ell,t})>\sigma(\varepsilon_{s,t})$), and seems preferable in all the parameters of our model. Nevertheless, the threshold inventory level for opening the long-distance market is still lower than that of the on-site market; the advanced demand information of the long-distance market is  still  valuable.

We have seen that the market preference of the firm can change as the current inventory position increases.  We study this phenomenon in Theorem \ref{lem_MultiplicativeDemand_samefuntion_IsIl} for the case of identical demand distributions and revenue functions in each market.

\begin{Theorem_}\label{lem_MultiplicativeDemand_samefuntion_IsIl}
If  $R_{s,t}(d)=R_{\ell,t}(d)$ for all $d$, $t$ and $\epsilon_{s,t}$, $\epsilon_{\ell,t}$ follow the same distribution, and  $\overline{\overline{d_s}}=\overline{\overline{d_\ell}}$, then following statements are true.
\begin{itemize}
  \item[i)] $d^*_{s,t}(I)\leq d^*_{\ell,t}(I)$ when  $ \E [ \epsilon_{s,t}H^{\prime}(I + q_t - \epsilon_{s,t} d^*_{s,t}(I))] \leq 0$.
  \item[ii))] $d^*_{s,t}(I)\geq d^*_{\ell,t}(I)$ when  $ \E [ \epsilon_{s,t}H^{\prime}(I + q_t - \epsilon_{s,t} d^*_{s,t}(I))] \geq 0$
\end{itemize}
\end{Theorem_}

Theorem~\ref{lem_MultiplicativeDemand_samefuntion_IsIl} establishes  the existence of a specific inventory level where if current inventory is below that level the demand intensity in the on-site market is optimally set higher than the demand intensity in the long-distance market.  When current inventory exceeds that specified level, then the long-distance market becomes the preferred sales channel.  The key insight is that inventory considerations have significant impact on the optimal selling strategy.  In limited inventory situations a retailer prefers selling through an online channel, whereas when inventory is plentiful the on-site market will be the preferred sales channel.  Due to the complexity of handling multiplicative demands, additional theoretical results seem intractable. However, the importance of this theorem is that it gives guidance to decision-making in a multi-market environment. The benefit of the long-distance market dominates the retailer's policy when the inventory level is low. While this paper only analyzes a one period potential for delayed shipment, if a longer delay for the long-distance market is possible, then the advantage of long-distance market is even greater (assuming low inventory).

\subsubsection{Correlated Demand}\label{SS_Perfect Correlated Demand}
In the previous sections, we made the assumption that demand in each market was independent. We now turn to the  significantly more complex case of correlated demands in the two markets.  For example, one would expect  that for certain  short lifecycle products, like fashion or high-tech products, demand in both markets will ebb and flow together with waning or rising consumer sentiment.
In the general case of correlated demands,  optimal policies may not posses a structure determined by  thresholds  and   insight for this  case is not easily achieved. However, we have obtained the result of  Proposition~\ref{lem_MultiplicativeDemand_Correlated_positivedst} below, that provides for a necessary condition for the on-site market to be opened.

\begin{Proposition_}\label{lem_MultiplicativeDemand_Correlated_positivedst}
When $D_{s,t}$ and $D_{\ell,t}$ are correlated, the following is true:  $$\mbox{If \ }  d_{s,t}^*(I)>0 , \mbox{ then  } V_t^{\prime}(I)<R_s^{\prime}(0).$$
\end{Proposition_}

Next we will discuss the perfect correlation case. Such cases are commonly used in studies of  inventory systems with multiple random variables, e.g,   \cite{Chen2013}.  For fixed $t$ the demands
$D_{s,t}$ and $ D_{\ell,t}$ exhibit perfect positive (negative) correlation when the variables
$\epsilon_{s,t}$ and $\epsilon_{\ell,t}$  are related with Eq. (\ref{eq:cor}),  for a constant  $a>0$ ($a<0$) .

\bea\label{eq:cor}
 \epsilon_{s,t}-1=a (\epsilon_{\ell,t}-1). \eea



In this case in Proposition~\ref{lem_MultiplicativeDemand_perfectcorrelated} below,  we  provide provable conditions regarding the on-site market's threshold inventory level for the case of perfectly correlated demand (negatively or positively correlated).

\begin{Proposition_}\label{lem_MultiplicativeDemand_perfectcorrelated}
When $d_{s,t}^*(I)=0$ and $\exists \delta>0$ satisfying that $d_{s,t}^*(I+\rho)>0$ for all $0<\rho<\delta$, the following are true:
\begin{itemize}
  \item $V_t^{\prime}(I)\leq R_s^{\prime}(0)$, if $D_{s,t}$ and $D_{\ell,t}$ exhibit perfect positive correlation;
  \item $V_t^{\prime}(I)\geq R_s^{\prime}(0)$,   if $D_{s,t}$ and $D_{\ell,t}$  exhibit perfect negative correlation.
\end{itemize}
\end{Proposition_}

Combining these two propositions, Theorem \ref{th_perfectnegativecorrelated_Chooseingon-siteStore_Threshold}
below, helps us shed light on comparing the correlated demand case and the   independent demand case. Specifically, if $D_{s,t}$ and $D_{\ell,t}$ exhibit perfect negative correlation, then the firm exhibits similar characterizations of the threshold inventory level for the on-site market.

\begin{Theorem_}\label{th_perfectnegativecorrelated_Chooseingon-siteStore_Threshold}
If $D_{s,t}$ and $D_{\ell,t}$ perfect negative correlated, then there exists a $I_{s,t}^*$ satisfying that $V_t^{\prime}(I_{s,t}^*)= R_s^{\prime}(0)$ and
\bea
d_{s,t}^*(I) =
\left\{
\begin{array}{cc}
>0 & I>I_{s,t}^*\\
=0 & I\leq I_{s,t}^*
\end{array}
\right.
\eea
\end{Theorem_}

\subsubsection{The Case in which Price is Constrained to be Equal in Both Markets}\label{SS_Unique_Price}

In this section, examine  the case where the two markets must charge the same price in each period as it is a policy many `clicks and mortar' retailers adopt. In this case the demand functions are specified by: $d_{\ell,t}=\beta_{\ell,t}\big(1-\frac{p}{\overline{\overline{p}}}\big)$ and $d_{s,t}=\beta_{s,t}\big(1-\frac{p}{\overline{\overline{p}}}\big)$. For period $t$, the optimization is written as follows:

\bea
V_t^u(I) &=& \max_{p \in [\underline{\underline{p}}, \overline{\overline{p}}]} J_t^u (I, p), \label{eqn_DP_Vut}
\eea
where
\bea
J_t^u (I, p) &=& \beta_{\ell,t}p\big(1-\frac{p}{\overline{\overline{p}}}\big) + \beta_{s,t}p\big(1-\frac{p}{\overline{\overline{p}}}\big) - \E H[I + q_t - \epsilon_{s,t}\beta_{s,t}\big(1-\frac{p}{\overline{\overline{p}}}\big)\nonumber\\
&&\, \, \, \, +  \alpha \E V_{t+1}[I + q_t - \epsilon_{s,t}\beta_{s,t}\big(1-\frac{p}{\overline{\overline{p}}}\big) -  \epsilon_{\ell,t}\beta_{\ell,t}\big(1-\frac{p}{\overline{\overline{p}}}\big)].
\label{eqn_DP_Jut}
\eea
Let $p^u_t(I)$ as the maximizer.

\begin{Lemma_}\label{lem_SS_Unique_Price}
$p^u_t(I)$ decreases in $I$.
\end{Lemma_}

Lemma~\ref{lem_SS_Unique_Price} proves the intuitive result that a retailer will decrease price in response to inventory level increases. Denote
$I_t^u*:=\max\{I|p^u_t(I)=\overline{\overline{p}}\}$. When the inventory level is lower than $I_t^u*$, the retailer will stop sales in both markets and hence, a threshold policy can be established.  Furthermore, we prove the following:

\begin{Theorem_}\label{th_SS_Unique_Sandwich}
$\min\{I_{s,t}^*, I_{\ell,t}^*\} \leq I_t^u*\leq  \max\{I_{s,t}^*, I_{\ell,t}^*\}$
\end{Theorem_}

Theorem~\ref{th_SS_Unique_Sandwich} places bounds on the threshold policy when a single price is used in both markets.  Interestingly, the retailer will open both markets at an inventory level in between the threshold inventory levels of the two markets in the previous price setting.  For example, sales will commence at an inventory level where the retailer in a heterogeneous price setting opens only one market.  Likewise, inventory levels where a retailer with greater pricing flexibility would close both markets will also represent an inventory level where the price constrained firm has closed both markets.

When a firm does have the ability to price differently in each market, the results of Theorem~\ref{th_SS_Unique_Sandwich} can be used to decrease the dimensionality of the retailer's pricing decision.  Insights, such as knowledge of $I_t^u*$, can provide a starting point for considering the market opening decision.  Given $I_t^u*$ and current inventory levels, the firm knows whether opening sales at both markets at the same price is optimal and hence, can use this information as a basis for considering which market should be opened earlier if able to charge different prices in the two markets.
 
 \begin{Remark_}\label{Remark_holding}
 
a) We note that when the holding cost approaches  zero, the difference  between a postponed  delivery and immediately shipping disappears. However, 
the analysis remains meaningful.  We note that following interesting special cases.

Case 1.  $\s=\h=0$ and $p_{st}(d) = p_{\ell t}(d)$ (i.e., $R_{st}(d)) =R_{\ell t}(d)) $,  then $I^*_{st}(d)= I^*_{\ell t} $ in Theorems 2 and    3 .

Case 2. $\s=0$,   $\h>0$ and $p_{st}(d) = p_{\ell t}(d)$ (i.e., $R_{st}(d)) =R_{\ell t}(d)) $,   
then $R_{s,t}^{\prime}(0) = R_{\ell,t}^{\prime}(0) $, 
 and if one goes through the steps of the proof of Theorem 3 one 
 can show that $I^*_{st} \le I^*_{\ell t} $. 
 Furthermore,  in this case if we assume the conditions of  
   Theorem 5 hold,  i.e.,  $\epsilon_{s,t}$, $\epsilon_{\ell,t}$ follow the same distribution,   then
  $ \E [ \epsilon_{s,t}H^{\prime}(I + q_t - \epsilon_{s,t} d^*_{s,t}(I))=\h>0$
and we have $d^*_{s,t}(I) \geq d^*_{\ell,t}(I).$ 

Case 3. $\s>0$,   $\h=0$ and $p_{st}(d) = p_{\ell t}(d)$ (i.e., $R_{st}(d)) =R_{\ell t}(d)) $,  then $R_{s,t}^{\prime}(0) = R_{\ell,t}^{\prime}(0) $, 
 and if one goes through the steps of the proof of Theorem 3 one 
 can show that $I^*_{st} \geq I^*_{\ell t} $. 
 In addition,  in this case if we assume the conditions of  
   Theorem 5 hold,  i.e.,  $\epsilon_{s,t}$, $\epsilon_{\ell,t}$ follow the same distribution,   then
  $ \E [ \epsilon_{s,t}H^{\prime}(I + q_t - \epsilon_{s,t} d^*_{s,t}(I))=-\s<0$
and we have $d^*_{s,t}(I) \leq d^*_{\ell,t}(I).$

b) Backlogging is adopted as a reasonable simplifying  assumption. 
We note that  for the primary market $s$    in our model (as in reality)     most of the time backlogging will 
   be incurred only for items demanded in the secondary market $\ell$. To see this consider the following realizations (examples), where we take   $T=2$ and $I_1=6$, $q_1=5$, $q_2=0$; recall that $I_{t+1}=I_t+q_t-D_{s,t}-D_{\ell,t}.$ 
  
 {\bf Example 1:}   $D_{s,1}=3$  $D_{\ell,1}=2$, $I_2=1+5$,  $D_{s,2}=2$  $D_{\ell,1}=2$.   Then   
 $I_3= 2 $,  i.e., no backlogged units at the end $T+1=3$. 
 
  {\bf Example 2:}   $D_{s,1}=5$  $D_{\ell,1}=3$, $I_2=-2+5=3$, $D_{s,2}=2$  $D_{\ell,1}=2$.  Then   
2 units of the $\ell$ market are  backlogged at the end of the 1st period and 
 $I_3= -1 $, i.e.,   $1$ unit of the $\ell$ market is  backlogged at the end $T+1=3$. 
 
 {\bf Example 3} (extreme and rare case):   $D_{s,1}=30 ,$  $D_{\ell,1}=2$, $I_2=-26+5=-21$,  $D_{s,2}=2$  $D_{\ell,1}=2$.  
 Then   
26 units (24 of the $s$ and 2  of the $\ell$ market) are  backlogged at the end of the 1st period and 
 $I_3= -20 $, i.e.,   $19$ units (from period 1)  are backlogged all the way to the end $(T+1=3)$ and 1 unit from period 2. 
Our computations confirm that for reasonable choices of the relevant costs, realizations  described in Example 3 do not occur.   However, allowing them   does simplify the analysis of this challenging problem. 
  
  \end{Remark_}

\section{Discussion and Conclusion}\label{S_Discussion2}
In this paper, we  investigated the structure of a retailer's optimal selling strategies when the retailer is faced with two distinct sales markets  and a common pool of replenishable  inventory with which to satisfy demand. A summary of the core selling strategies derived in \S\ref{S_Discussion} is provided in Table~\ref{Table_Summary}.  The optimal policies in this  table indicate that  both  inventory and demand uncertainty are key issues   in choosing which markets to open up for demand, i.e., choose a price for which the mean market demand is positive.
Based on the type of demand uncertainty, it is shown that optimal policies are characterized by inventory level thresholds
which determine whether which market  will  offer  the product  for sale and at what price.

The following managerial  insight is provided.  When demand uncertainty in both markets is additive and inventory is low, selling exclusively through the online market is preferred. As inventory increases beyond a specific threshold, the online market is also opened (offers the item for sale), and the optimal selling strategy is to increase demand intensity in both markets as $I$ increases further. However, the problem is more  challenging  when both markets have  multiplicative demand uncertainty. In this case it  is shown that the sales in one channel may actually decrease as the inventory level increases. Despite this, we show  that the demand intensity in at least one of the markets  does indeed increase as the inventory level increases.  A key result is provided by the following property we establish in the multiplicative mean demand model case. In this case it is shown that optimal  threshold policies can derived from benchmark problems where demand uncertainty in the current period can be ignored.  Hence, if a retailer would sell through a particular market assuming no uncertainty in demand in the current period, then they should also open that market when multiplicative uncertainty exists.

\begin{table}[ht]
\caption{Optimal Selling Strategies under Different Demand Uncertainty Models}\label{Table_Summary}
\begin{center}
\resizebox{1.0\textwidth}{!}{
\begin{tabular}{||c|c||l||}
\hline
\multicolumn{2}{||c||}{Type of Demand Certainty } & \\
\cline{1-2}
on-site market & long-distance market & \multicolumn{1}{|c||}{Optimal selling strategy}\\
\hline
& & There exist an inventory level $\bar I_t$ and a pair of mean demands $(\bar d_{s,t},\bar d_{\ell,t})$.\\
additive& additive& If $I\geq \bar I_t$, then $d_{s,t}^*(I)$, $\bar I_t-d_{s,t}^*(I)$ and $d_{\ell,t}^*(I)$ increasing in $I^*$.\\
& & If $I\leq \bar I_t$, then $d_{s,t}^*(I)=0$ and $d_{\ell,t}^*(I)$ increases as $I$ increases.\\
\hline
& & There exist an inventory level $\bar I_t$ and a pair of the mean demands $(\bar d_{s,t},\bar d_{\ell,t})$.\\
additive& multiplicative& If $I\geq \bar I_t$, then $d_{\ell,t}^*(I)$, $\bar I_t-d_{\ell,t}^*(I)$ and $d_{s,t}^*(I)$ increasing in $I^*(I)$.\\
& & If $I\leq \bar I_t$, then $d_{\ell,t}^*(I)=0$ and $d_{s,t}^*(I)$ increases as $I$ increases.\\
\hline
& & Optimal selling policy is a threshold policy. \\
multiplicative& additive& Optimal selling quantities at both on-site market and long-distance market increase as the inventory level increases.\\
\hline
& & Optimal selling policy is a threshold policy. \\
multiplicative& multiplicative& Optimal selling quantity at on-site market or long-distance market may increase or decrease as the inventory level increases. \\
& & Optimal selling quantity at at least one of on-site market and long-distance market increases as inventory level increases.\\
\hline
\end{tabular}}
\end{center}
\end{table}

For a retailer serving two markets from a common pool of inventory, the insights from this paper can be valuable.
With limited inventory, the retailer will prefer to sell through only one market.  In this limited inventory scenario, analyzing cases where the marginal profit in each market is equal reveals an often overlooked benefit of selling through the online channel; namely, gaining  demand information is a truly valuable aspect of serving this market.  Evidence of this was shown through characterization of cases where a firm prefers sales through the long-distance market despite the marginal profit of the on-site market being significantly larger.


 In conclusion, this study  contributes to the literature by providing 
 the first model for a firm's    dual market dynamic   pricing  problem 
 in the presence of  exogenously determined  inventory replenishment considerations.       
      Our study of a finite horizon case is consistent with fashion products, whereas the infinite horizon case might provide tractable insight for functional products.

In current and future work, we plan to investigate the  effect of relaxing   some 
of  our modeling assumptions.  For instance, a retailer may choose not to handle online demand in a batch at the end of the period and exploration of flexibility in the long-distance market's timing of demand satisfaction (e.g. satisfy in this period or the subsequent one) would provide meaningful results.    Extensions treating the lost sales case and providing more thorough treatment of the correlation of the two market's demands  are also  the subject of future research.    Herein, we  have provided  conditions regarding the on-site market's threshold inventory level for the case of perfectly  (negatively or positively)  correlated demands.     In the general case, if the two markets demands are correlated, then the threshold structure of an optimal  policy may not  hold.

\appendix
\newtheorem{LemmaA_}{Lemma}[section]
\newtheorem{DefinitionA_}{Definition}[section]

\section{Appendix: Proofs}\label{SA_Proofs}

\newcommand{\RL}{{\rm I\hspace{-0.8mm}R}}

%
%
In this section we will use the following   simplified definition of supermodularity (i.e., ``increasing differences'') of a real function on $\RL^2$  given  in Definition~\ref{Def_Supermodular} below. 
A comprehensive exposition of the topic of supermodularity and its application in sequential decision problems is given in \cite{Topkis1998}.
\begin{DefinitionA_}\label{Def_Supermodular}
A real function $f(x,y)$ is supermodular in $(x,y) \in \RL^2$ if
  $$f(x_1,y_1)+f(x_2,y_2)\geq f(x_1,y_2)+f(x_2,y_1)$$ for all $x_1\geq x_2$ and $y_1\geq y_2$.
 \end{DefinitionA_}
If the inequality above is reversed, the function $f$ is called \textit{submodular.}

Lemma~\ref{lem_Supermodular} below summarizes properties of supermodular functions that we will use, we refer the reader to \cite{Topkis1998} (Theorem 2.6.2) and \cite{simchi2005logic} for proofs.
\begin{LemmaA_}\label{lem_Supermodular}
{(Topkis 1998)}
\begin{enumerate}
\item If $g(x)$ $\RL \rightarrow \RL $ is concave, then $$f(x,y)=g(x+y)$$ is submodular in $(x,y)$ and  $$h(x,y)=g(x-y)$$ is supermodular in $(x,y)$.
 \item If $f(x,y)$ is supermodular (submodular) in $(x,y) \in \RL^2$ and $$\tilde x^*_{f}(y)=arg\,max_{x} f(x,y)$$  then $\tilde x^*_{f}(y)$ is increasing (decreasing) in $y$.
\end{enumerate}
\end{LemmaA_}


\begin{pf}{\bf of Lemma~\ref{lem_ConcaveRevenuefunction}.} We first clarify the existence of optimality. The existence of the optimality follows  from the continuity of the value function and boundary condition of $d\in[\underline{\underline{d}},\overline{\overline{d}}]$. Denote $c^H=\max\{\s,\h\}$.
First, we have $V_{T+1}(I)=c_{e}I^{+}>0.$ Also,
$\E V_{T+1}(I_{T}+q_{T}-D_{s,T}-D_{\ell,T})=c_{e}\E(I_{T}+q_{T}-D_{s,T}-D_{\ell,T})^{+}\ge$ $c_{e}(I_{T}+q_{T}- (d_{s,T}+d_{\ell,T}))>-\infty .$
The proof is easy to complete using induction and Eq. (\ref{eqn_DP_Vt}).

Second, we will show that $V_t$ is concave in $I$. From the terminal condition that $V_{T+1}(I)=\e \min \{I,0\}$ for all $I$, it is clear $V_{T+1}$ is concave. If we assume $V_t$ is concave, then because  the linear combination of two concave functions is concave, we have that $J_t(I, d_s, d_\ell)$ is concave in  $(I, d_s, d_\ell)$. Given that concavity is preserved under maximization \citep{HeymanSobel1984}, we obtain that $V_{t-1}$ is concave. The second statement of Lemma~\ref{lem_ConcaveRevenuefunction} also follows by induction.  The last statement follows from 
the concavity of $V$, $-H$, $R_s$ and $R_\ell$ and Lemma ~\ref{lem_Supermodular}.
\end{pf}

\begin{pf}{\bf of Lemma~\ref{lem_GeneralCase_Increasing}.} We only need to show that if $\delta>0$ and $d_{s,t}^*(I+\delta) < d_{s,t}^*(I)$ then $d_{\ell,t}^*(I+\delta) \geq d_{\ell,t}^*(I)$. We can write
%
$d_{\ell,t}^*(I)=\displaystyle arg\,\!\!\max_{\substack{\underline{\underline{d_\ell}}\leq d_\ell\leq \overline{\overline{d_\ell}}}} J_t(I, d_{s,t}^*(I), d_\ell)$ as\\
$$\ = \displaystyle arg\,\!\!\max_{\substack{\underline{\underline{d_\ell}}\leq d_\ell\leq \overline{\overline{d_\ell}}}} R_s(d_{s,t}^*(I)) + R_{\ell,t}(d_\ell) - \E H(I+q_t-\epsilon_{s,t} d_{s,t}^*(I)) + \alpha \E V_{t+1}(I+q_t- \epsilon_{s,t} d_{s,t}^*(I) -  \epsilon_{\ell,t} d_\ell).
 $$

Because $-H$ and $V_{t+1}$ are concave, $J_t(I, d_{s,t}^*(I), d_\ell)$ is supermodular in $(I,d_\ell)$ and submodular in $(d_{s,t}^*(I), d_\ell)$. Now,  from Lemma A.1, we have $d_\ell^*(I)$ increases in $I$ and decreases in $d_{s,t}^*$. Hence, from $I + \delta > I$ and $d_{s,t}^*(I+\delta) < d_{s,t}^*(I)$, we can obtain $d_{\ell,t}^*( I + \delta) \geq d_{\ell,t}^*(I)$.
\end{pf}

%
%
\begin{pf} {\bf of Theorem~\ref{lem_Additiveon-siteDemand}.}
\emph{To see part $i)$},  from (\ref{eqn_DP_Vt}) and (\ref{eqn_DP_Jt}), we first have that
\bea
V_t(I) &=& \max\limits_{\substack{\underline{\underline{d_s}}\leq d_s\leq \overline{\overline{d_s}}}}
 \big[R_{s,t}(d_s)  - \E H(I + q_t -  d_s - \omega_{s,t})
 + f(I - d_s)\label{eqn_lem_AdditiveDemandStore_q1}
 \big]
\eea
\bea
\mbox{where } &&f(x)=\max\limits_{\substack{\underline{\underline{d_\ell}}\leq d_\ell\leq \overline{\overline{d_\ell}}}}
 \big[
 R_{\ell,t}(d_\ell) + \alpha \E V_{t+1}(x + q_t- \omega_{s,t} - \epsilon_{\ell,t} d_\ell - \omega_{\ell,t})
 \big]\label{eqn_lem_AdditiveDemandStore_q2}
\eea
Because $V_t$ is concave and concavity is preserved under maximization, $f$ is concave. We can show that the right side of Eq. (\ref{eqn_lem_AdditiveDemandStore_q1}) is supermodular in  $(I,d_s)$. The first term depends only on   $d_s$ and it is clearly submodular in $(I,d_s)$. The supermodularity of the second and the third terms follows from the concavity of $-H$, $f$ and linear combination of $I-d_s$
Therefore, from from Lemman \ref{lem_Supermodular}  we have $d_{s,t}^*(I)$ increases in $I .$

\emph{To see part ii)}, we denote $x=I - d_s$ and $x^*(I)=I-d_{s,t}^*(I)$. From (\ref{eqn_DP_Vt}), we have
\bea
V_t(I) &=& \max\limits_{\substack{I+\overline{\overline{d_s}}\leq x\leq I+ \underline{\underline{d_s}}}}
 \big[R_{s,t}(x-I)  - \E H(x + q_t - \omega_{s,t})
 + f(x) \big] \label{eqn_lem_AdditiveDemandStore_q3}
\eea
Following  a similar argument as in the proof of  part $i)$, we can show that   $x^*(I)$ increases in $I .$

To establish  part iii), we note  that the right side of Eq.  (\ref{eqn_lem_AdditiveDemandStore_q2}) is supermodular in $(x,d_\ell)$. The first term depends only on $d_\ell$ and it is clearly supermodular in $(x,d_\ell)$. The supermodularity of the second term follows from the concavity of $V_t$ and linear combination of $x-\epsilon_{\ell,t} d_\ell$. Therefore, we have $d_{\ell,t}^*(I)$ is increasing  in $x^*(I)$. Combining this with the result of part $ii)$, we have that  $d_{\ell,t}^*(I)$ increases in $I$.
\end{pf}

\begin{pf}{\bf of Remark~\ref{Remark_Additiveon-siteDemand_monotonIn_qt}.}
 At time $t$ consider two forthcoming order quantities  $q_t^a $, $ q_t^b$  with  $q_t^a > q_t^b$ and let  $V^j$, $J_t^j$ and $(d_{s,t}^{*,j}(I),d_{\ell,t}^{*,j}(I))$ denote the optimal profit, objective function and optimal selling price corresponding to the  forthcoming order quantities $q_t^j$, $j=a, b$. From Eq. (\ref{eqn_DP_Jt}) it follows that $J_t^a(I,d_s,d_\ell)=J_t^b(I + q_t^a -q_t^b,d_s,d_\ell)$ for all $d_s$ and $d_\ell$.  Now using Lemma ~\ref{lem_Additiveon-siteDemand},  we obtain: $d_{i,t}^{*,a}(I) = d_{i,t}^{*,b}(I + q_t^a -q_t^b)\geq d_{i,t}^{*,b}(I)$ for $i=s,\ell$.
\end{pf}

\begin{pf}{\bf of Lemma ~\ref{lem_BL_BS_Monotone}.}
To prove part i), we use the first order partial derivatives. From the supermodularity of $J_t^{B_\ell*}$ in $(I,d_\ell)$, we have $d_\ell^{B_\ell*}(I)$ increases in $I$ if $d_s\in \{0, \overline{\overline{d_s}}\}$. For the case in which
$d_s\in (0, \overline{\overline{d_s}})$, we will show that $\frac{d d_\ell^{B_\ell*}(I)}{dI}\geq 0$. Because $d_s>0$ and $\epsilon_{s,t}$ is a random variable with  a continuous distribution, $I + q_t - \epsilon_{s,t}d_s\neq 0$ almost everywhere, we only need to focus on the interval where $H'$ and $H''$ exist and we note that $H''(I + q_t - \epsilon_{s,t}d_s)=0$ almost everywhere.

From the first order partial derivatives, we have that when $d_s=d_{s,t}^{B_\ell*}(I)$ and $d_\ell=d_\ell^{B_\ell*}(I)$,
\bean
&&0=R_{s,t}^{\prime}(d_s)  +  \E\big[ \epsilon_{s,t} H^{\prime}( I + q_t - \epsilon_{s,t}d_s)\big] - \alpha \E \big[\epsilon_{s,t} V_{t+1}^{\prime}(I + q_t - \epsilon_{s,t} d_s - d_\ell )\big]
 \label{eqn_eqn_lm_BL_BS_Monotone_q02}
\\
&&0=R_{\ell,t}^{\prime}(d_\ell) - \alpha \E \big[V_{t+1}^{\prime}(I + q_t - \epsilon_{s,t} d_s-d_\ell)\big] \label{eqn_eqn_lm_BL_BS_Monotone_q03}
\eean
and with $H^{\prime \prime}(I + q_t - \epsilon_{s,t}d_s)=0$ almost everywhere, we have
\bea
&0&=R_{s,t}^{\prime\prime}(d_{s,t}^{B_\ell*}(I))\frac{d d_{s,t}^{B_\ell*}(I)}{dI}
- \alpha \E \big[\epsilon_{s,t}(1 - \epsilon_{s,t}\frac{d d_{s,t}^{B_\ell*}(I)}{dI} - \frac{d d_\ell^{B_\ell*}(I)}{dI}) V_{t+1}^{\prime\prime}(I + q_t - \epsilon_{s,t} d_{s,t}^{B_\ell*}(I) -  d_\ell^{B_\ell*}(I))\big]\nonumber \\
 \label{eqn_eqn_lm_BL_BS_Monotone_q04}
\eea
\bea
&0&= R_{\ell,t}^{\prime\prime}(d_\ell^{B_\ell*}(I))\frac{d d_\ell^{B_\ell*}(I)}{dI} - \alpha \E \big[(1 - \epsilon_{s,t}\frac{d d_{s,t}^{B_\ell*}(I)}{dI}
- \frac{d d_\ell^{B_\ell*}(I)}{dI})V_{t+1}^{\prime\prime}(I + q_t - \epsilon_{s,t} d_{s,t}^{B_\ell*}(I) -  d_\ell^{B_\ell*}(I))\big]\nonumber\\
\label{eqn_eqn_lm_BL_BS_Monotone_q05}
\eea
To simplify the notation, denote
\bean
&d_s^{\prime}(I)=\frac{d d_{s,t}^{B_\ell*}(I)}{dI},  \,\,\,\,\,\,\,\,
d_\ell^{\prime}(I)=\frac{d d_\ell^{B_\ell*}(I)}{dI},\\
&\gamma_s = -R_{s,t}^{\prime\prime}(d_{s,t}^{B_\ell*}(I)),\,\,\,\,\,\,\,\,
\gamma_\ell =-R_{\ell,t}^{\prime\prime}(d_\ell^{B_\ell*}(I)),\\
&\varpi_{\epsilon_{s,t}} =-V_{t+1}^{\prime\prime}(I + q_t - \epsilon_{s,t} d_{s,t}^{B_\ell*}(I) -  d_\ell^{B_\ell*}(I))\,\,\,\,\,  \forall \epsilon_{s,t}.
\eean
From the strict concavity of $R_{s,t}$ and  $R_{\ell,t}$ and the concavity of $-H$ and $V$, we have
\bea
\gamma_s >0, \,\,\,\gamma_\ell >0,  \,\,\,\varpi_{\epsilon_{s,t}} \geq 0.
\label{eqn_eqn_lm_BL_BS_Monotone_q06}
\eea
We can rewrite  (\ref{eqn_eqn_lm_BL_BS_Monotone_q04}) and (\ref{eqn_eqn_lm_BL_BS_Monotone_q05}) as follows.
\bean
&&\big[\gamma_s +\alpha \E(\epsilon_{s,t}^2\varpi_{\epsilon_{s,t}} )\big] d_s^{\prime}(I) + \alpha \E(\epsilon_{s,t}\varpi_{\epsilon_{s,t}} ) d_\ell^{\prime}(I) = \alpha \E(\epsilon_{s,t}\varpi_{\epsilon_{s,t}} ) \label{eqn_eqn_lm_BL_BS_Monotone_q07}\\
&&\alpha \E(\epsilon_{s,t}\varpi_{\epsilon_{s,t}} ) d_s^{\prime}(I) + \big[\gamma_\ell+\alpha \E(\varpi_{\epsilon_{s,t}} )\big] d_\ell^{\prime}(I) = \alpha \E(\varpi_{\epsilon_{s,t}} ) \label{eqn_eqn_lm_BL_BS_Monotone_q08}
\eean
From above two equations, we have
\bea
d_\ell^{\prime}(I)=
\frac{  \alpha \E(\varpi_{\epsilon_{s,t}} )\big[\gamma_s  +\alpha \E(\epsilon_{s,t}^2\varpi_{\epsilon_{s,t}} )\big] - \alpha \E(\epsilon_{s,t}\varpi_{\epsilon_{s,t}} )\alpha \E(\epsilon_{s,t}\varpi_{\epsilon_{s,t}} ) }
{\big[\gamma_\ell+\alpha \E(\varpi_{\epsilon_{s,t}} )\big]\big[\gamma_s  +\alpha \E(\epsilon_{s,t}^2\varpi_{\epsilon_{s,t}} )\big]-\alpha^2 \E(\epsilon_{s,t}\varpi_{\epsilon_{s,t}} )^2}
\label{eqn_eqn_lm_BL_BS_Monotone_q09}
\eea
From the Cauchy-Schwarz inequality, we have
\bea
\E(\varpi_{\epsilon_{s,t}} )\E(\epsilon_{s,t}^2\varpi_{\epsilon_{s,t}} )- \E^2(\epsilon_{s,t}\varpi_{\epsilon_{s,t}})
= \E(\sqrt{\varpi_{\epsilon_{s,t}}})^2 \E(\epsilon_{s,t}\sqrt{\varpi_{\epsilon_{s,t}}} )^2 - \E^2[(\epsilon_{s,t}\sqrt{\varpi_{\epsilon_{s,t}}})\sqrt{\varpi_{\epsilon_{s,t}}} ]
\geq 0.
\label{eqn_eqn_lm_BL_BS_Monotone_q10}
\eea
Combining this result with Ineqs. (\ref{eqn_eqn_lm_BL_BS_Monotone_q06}), we have:
%
%
\begin{align*}
\ & \alpha \E(\varpi_{\epsilon_{s,t}} )\big[\gamma_s  + \alpha \E(\epsilon_{s,t}^2\varpi_{\epsilon_{s,t}} )\big] - \alpha \E(\epsilon_{s,t}\varpi_{\epsilon_{s,t}} )\alpha \E(\epsilon_{s,t}\varpi_{\epsilon_{s,t}} )  \ & \ \\
 & \geq
\alpha^2\big[\E(\varpi_{\epsilon_{s,t}} )\E(\epsilon_{s,t}^2\varpi_{\epsilon_{s,t}} )- \E^2(\epsilon_{s,t}\varpi_{\epsilon_{s,t}} )\big] \geq 0, &\ & &
\end{align*}
and
\bean
\big[\gamma_\ell+\alpha \E(\varpi_{\epsilon_{s,t}} )\big]\big[\gamma_s  +\alpha \E(\epsilon_{s,t}^2\varpi_{\epsilon_{s,t}} )\big]-\alpha^2 \E(\epsilon_{s,t}\varpi_{\epsilon_{s,t}} )^2
>\alpha^2\big[\E(\varpi_{\epsilon_{s,t}} )\E(\epsilon_{s,t}^2\varpi_{\epsilon_{s,t}} )
\geq 0 .
\eean
The above together with Eq. (\ref{eqn_eqn_lm_BL_BS_Monotone_q09}), imply that  $d_\ell^{\prime}(I)\geq 0$. That is to say, $\frac{d d_\ell^{B_\ell*}(I)}{dI}\geq 0$.

To see Part (ii), we rewrite objective function as follows
\bean
V_t^{B_s}(I)=\max_{d_s} R_{s,t}(d_s) + G(I-d_s)
\eean
where
\bean
G(x)=\max_{d_\ell} R_{\ell,t}(d_\ell) -\E H(x + q_t) + \alpha \E V_{t+1}(x + q_t  - \epsilon_{\ell,t}d_\ell).
\eean
First, we have $G(x)$ is concave in $x$. Hence, $G(I-d_s)$ is supermodular in $(I,d_s)$. Further, \\$R_{s,t}(d_s) + G(I-d_s)$ is supermodular in $(I,d_s)$. Therefore, $d_s^{B_s*}(I)$ is increasing in $I$.
\end{pf}

\begin{pf}{\bf of Lemma ~\ref{Lem_BL_BS_dlBLanddl}.} First, we show that if $d_{\ell,t}^{B_\ell*}(I)=0$ then $d_{\ell,t}^*(I)=0$. For all $d_s$, $d_\ell$, we have
\bea
J_t^{B_\ell}(I,d_s,d_\ell)-\! J_t(I,d_s,d_\ell)
\!=\! \alpha \E V_{t+1}( I + q_t \!- \!\epsilon_{s,t} d_s- \!d_\ell) \!-\!\alpha
 \E V_{t+1}( I + q_t \! - \!\epsilon_{s,t}d_s- \!\epsilon_{s,t}d_\ell) \geq 0. \label{eqn_lm_BL_BS_dlBLanddl_q1}
\eea
The inequality follows from the concavity of $V_{t+1}$ and Jensen's inequality.

If $d_{\ell,t}^{B_\ell*}(I)=0$, we can obtain that
\bea
J_t(I,d_{s,t}^{B_\ell*}(I),0) = J_t^{B_\ell}(I,d_{s,t}^{B_\ell*}(I),0)\geq J_t^{B_\ell}(I,d_s,d_\ell)\geq J_t(I,d_s,d_\ell) \,\,\,\,\, \forall d_s, d_\ell .
\label{eqn_lm_BL_BS_dlBLanddl_q2}
\eea
The first inequality follows from the optimality of $(d_{s,t}^{B_\ell*}(I),0)$. The second inequality follows from Ineq. (\ref{eqn_lm_BL_BS_dlBLanddl_q1}). Notice that we always pick up the the maximizer with the smallest $d_\ell^*$. Therefore, we have  $d_{\ell,t}^*(I)=0$ from (\ref{eqn_lm_BL_BS_dlBLanddl_q2}).

Below, we will show that if $d_{\ell,t}^*(I)=0$, then $d_{\ell,t}^{B_\ell*}(I)=0$ and prove the statement by contradiction. Suppose $d_{\ell,t}^*(I^a)=0$ and $d_{\ell,t}^{B_\ell*}(I^a)>0$. From $d_{\ell,t}^*(I)=0$, we can find $I^b\geq I^a$ and a positive number $\delta>0$ satisfying the conditions that $d_{\ell_t}^{*}(I^b)=0$ and $d_{\ell_t}^{*}(I^b+\rho)>0$ for all $0<\rho<\delta$. We still notice that

\begin{align}
\frac{\partial J_t}{\partial d_s} |_{d_\ell=0}
&=R_{s,t}^{\prime}(d_s)  +  \E [ \epsilon_{s,t} H^{\prime}( I + q_t - \epsilon_{s,t}d_s) ] - \alpha \E  [\epsilon_{s,t} V_{t+1}^{\prime}(I + q_t - \epsilon_{s,t} d_s  ) ] \nonumber  \\
&=\frac{\partial J_t^{B_\ell}}{\partial d_s}|_{d_\ell=0}\label{eqn_lm_BL_BS_dlBLanddl_q3}
\end{align}

\begin{align}
\frac{\partial J_t}{\partial d_\ell} |_{d_\ell=0}
&=R_{\ell,t}^{\prime}(0) - \alpha \E [ \epsilon_{\ell,t} V_{t+1}^{\prime}(I + q_t - \epsilon_{s,t} d_s)]
=R_{\ell,t}^{\prime}(0) - \alpha \E  [ V_{t+1}^{\prime}(I + q_t - \epsilon_{s,t} d_s) ] \nonumber  \\
&=\frac{\partial J_t^{B_\ell}}{\partial d_\ell}|_{d_\ell=0}\label{eqn_lm_BL_BS_dlBLanddl_q4}
\end{align}
%
%
%
Hence, we have
\bea
\frac{\partial J_t^{B_\ell}(I^b,d_s,d_\ell)}{\partial d_\ell}\big|_{(d_s,d_\ell)=(d_s^*(I^b),0)}=\frac{\partial J_t(I^b,d_s,d_\ell)}{\partial d_\ell}\big|_{(d_s,d_\ell)=(d_s^*(I^b),0)} =0 \label{eqn_lm_BL_BS_dlBLanddl_q5}
\eea

 If $d_{s,t}^*(I^b)=0$, then combining with (\ref{eqn_lm_BL_BS_dlBLanddl_q3}) we have
$$
\frac{\partial J_t^{B_\ell}(I^b,d_s,d_\ell)}{\partial d_s}\big|_{(d_s,d_\ell)=(0,0)} =\frac{\partial J_t(I^b,d_s,d_\ell)}{\partial d_s}\big|_{(d_s,d_\ell)=(0,0)}\leq 0.$$ Together with Eq. (\ref{eqn_lm_BL_BS_dlBLanddl_q5}), we have  $d_{\ell,t}^{B_\ell*}(I^b)=0$ and  $d_{s,t}^{B_\ell*}(I^b)=0$.

If $d_{s,t}^*(I^b)=\overline{\overline{d_s}}$, then combining with Eq. (\ref{eqn_lm_BL_BS_dlBLanddl_q3}) we will have
$$
\frac{\partial J_t^{B_\ell}(I^b,d_s,d_\ell)}{\partial d_s}\big|_{(d_s,d_\ell)=(\overline{\overline{d_s}},0)} =\frac{\partial J_t(I^b,d_s,d_\ell)}{\partial d_s}\big|_{(d_s,d_\ell)=(\overline{\overline{d_s}},0)}\geq 0
.$$  Together with Eq. (\ref{eqn_lm_BL_BS_dlBLanddl_q5}),   we have $d_{\ell,t}^{B_\ell*}(I^b)=0$ and  $d_{s,t}^{B_\ell*}(I^b)=\overline{\overline{d_s}}$.

If $d_{s,t}^*(I^b)\in (0, \overline{\overline{d_s}})$, then  combining with (\ref{eqn_lm_BL_BS_dlBLanddl_q3}) we have
$$\frac{\partial J_t^{B_\ell}(I^b,d_s,d_\ell)}{\partial d_s}\big|_{d_s=d_{s,t}^*(I^b), d_\ell=0} =\frac{\partial J_t(I^b,d_s,d_\ell)}{\partial d_s}\big|_{ d_s=d_{s,t}^*(I^b), d_\ell=0}=0.$$  Together with Eq. (\ref{eqn_lm_BL_BS_dlBLanddl_q5}),   we have $d_{\ell,t}^{B_\ell*}(I^b)=0$ and  $d_{s,t}^{B_\ell*}(I^b)=d_{s,t}^*(I^b)$.

In the above three cases, we can get $d_{\ell,t}^{B_\ell*}(I^b)=0$. However, from Lemma ~\ref{lem_BL_BS_Monotone}, we have $d_{\ell,t}^{B_\ell*}(I^b)\geq d_{\ell,t}^{B_\ell*}(I^a)>0$ and thus, have a contradiction. Therefore, we reach $d_{\ell,t}^{B_\ell*}(I)=0$ if and only if $d_{\ell,t}^*(I)=0$.

To see part (ii), we can get that $d_{s,t}^{B_s*}(I)=0$ if and only if $d_{s,t}^*(I)=0$ by similar argument as the part (i).
\end{pf}

\begin{pf} {\bf of Theorem ~\ref{th_MultiplicativeDemand_Threshold}.}
From Lemmas \ref{lem_BL_BS_Monotone} and \ref{Lem_BL_BS_dlBLanddl}, there exist $I_{\ell,t}^*$ satisfying that $d_{\ell,t}^{B_\ell*}(I)>0$, $d_{\ell,t}^{*}(I)>0$ if $I>I_{\ell,t}^*$ and $d_{\ell,t}^{*}(I)=d_{\ell,t}^{B_\ell*}(I)=0$ if $I\leq I_{\ell,t}^*$. Follow the same logic, we have
there exist $I_{s,t}^*$ satisfying that $d_{s,t}^{B_s*}(I)>0$, $d_{s,t}^{*}(I)>0$ if $I>I_{s,t}^*$ and $d_{s,t}^{*}(I)=d_{s,t}^{B_s*}(I)=0$ if $I\leq I_{s,t}^*$.

From the envelop theorem and (
5), we also have that when $I=I_{s,t}^*$,
\bean
&&V^{\prime}(I) = H^{\prime}(I + q_t) + \alpha \E \big[ V_{t+1}^{\prime}(I + q_t -  \epsilon_{\ell,t} d_{\ell}^*(I))\big]\\
&&0=\frac{\partial J_t}{\partial d_s} = R_{s,t}^{\prime}(0) -  H^{\prime}(I + q_t) - \alpha \E \big[ V_{t+1}^{\prime}(I + q_t -  \epsilon_{\ell,t} d_{\ell}^*(I))\big]
\eean
Hence,  $V_t^{\prime}(I_{s,t}^*)= R_{s,t}^{\prime}(0)$.
\end{pf}

\begin{pf}{\bf of Lemma ~\ref{lem_MultiplicativeDemand_ChooseingOnline_Threshold_Vprime}.}
From $d_{\ell,t}^*(I_{\ell,t}^*)=0$, we have
\bea  \label{eqn_lm_MultiplicativeDemand_ChooseingOnline_Threshold_Vprime_q1}
V_t^{\prime}(I_{\ell,t}^*)= -  \E \big[ H^{\prime}(I_{\ell,t}^* + q_t - \epsilon_{s,t} d_{s,t}^*(I_{\ell,t}^*))\big] + \alpha \E \big[ V_{t+1}^{\prime}(I_{\ell,t}^* + q_t - \epsilon_{\ell,t} d_{\ell,t}^*(I_{\ell,t}^*))\big].
\eea

From the  first order partial derivatives and Theorem ~\ref{th_MultiplicativeDemand_Threshold}, we also have
\bea
0&&=\frac{\partial J_t(I,d_s,d_\ell)}{\partial d_\ell}|_{d_\ell=0}\nonumber\\
&&=R_\ell^{\prime}(0)  - \alpha \E \big[\epsilon_{\ell,t} V_{t+1}^{\prime}(I_{\ell,t}^* + q_t -\epsilon_{s,t} d_{s,t}^*(I_{\ell,t}^*) - \epsilon_{\ell,t} d_{\ell,t}^*(I_{\ell,t}^*))\big]  \nonumber\\
&&= R_\ell^{\prime}(0)  - \alpha \E V_{t+1}^{\prime}(I_{\ell,t}^* + q_t -\epsilon_{s,t} d_{s,t}^*(I_{\ell,t}^*)) .\label{eqn_lm_MultiplicativeDemand_ChooseingOnline_Threshold_Vprime_q2}
\eea
Combining Eqs. (\ref{eqn_lm_MultiplicativeDemand_ChooseingOnline_Threshold_Vprime_q1}) and (\ref{eqn_lm_MultiplicativeDemand_ChooseingOnline_Threshold_Vprime_q2}), we have $V_t^{\prime}(I_{\ell,t}^*)= R_\ell^{\prime}(0) -  \E \big[ H^{\prime}(I_{\ell,t}^* + q_t - \epsilon_{s,t} d_{s,t}^*(I_{\ell,t}^*))\big]$.
\end{pf}

\begin{pf}{\bf of Theorem ~\ref{lem_MultiplicativeDemand_IsIl}.} If $R_{s,t}^{\prime}(0) > R_{\ell,t}^{\prime}(0) + c_p$, then using  Theorem~\ref{th_MultiplicativeDemand_Threshold}  and Lemma ~\ref{lem_MultiplicativeDemand_ChooseingOnline_Threshold_Vprime} we have that $V_t^{\prime}(I_{s,t}^*)=R_{s,t}^{\prime}(0) > R_{\ell,t}^{\prime}(0) + c_p \geq R_{\ell,t}^{\prime}(0) -  \E \big[ H^{\prime}(I_{\ell,t}^* + q_t - \epsilon_{s,t} d_{s,t}^*(I_{\ell,t}^*))\big] = V_t^{\prime}(I_{\ell,t}^*)$. Since $V_t$ is concave, $I_{s,t}^*<I_{\ell,t}^*$.

If $R_{\ell,t}^{\prime}(0) > R_{s,t}^{\prime}(0) + h$, then using Theorem ~\ref{th_MultiplicativeDemand_Threshold}  and Lemma~\ref{lem_MultiplicativeDemand_ChooseingOnline_Threshold_Vprime} we obtain  $V_t^{\prime}(I_{s,t}^*)=R_{s,t}^{\prime}(0) < R_{\ell,t}^{\prime}(0) - c_h \leq R_{\ell,t}^{\prime}(0) -  \E \big[ H^{\prime}(I_{\ell,t}^* + q_t - \epsilon_{s,t} d_{s,t}^*(I_{\ell,t}^*))\big] = V_t^{\prime}(I_{\ell,t}^*)$. Because $V_t$ is concave, we have  $I_{s,t}^*>I_{\ell,t}^*$.
\end{pf}

\begin{pf}{\bf of Theorem ~\ref{lem_MultiplicativeDemand_preferIl}.} The first statement follows directly from Theorem ~\ref{lem_MultiplicativeDemand_IsIl} and thus, we focus on proving the second statement. Instead of proving the second statement directly, we will use mathematic induction to prove that if $R_{s,t}^{\prime}(0) < R_{\ell,t}^{\prime}(0) + \s$, then for all $1\leq t\leq T-1$,
\bea
I_{\ell,t}^*<I_{s,t}^* \mbox{ and } V_t^{\prime}(0)\leq R_{\ell,t}^{\prime}(0) + \s .
\label{eqn_lm_MultiplicativeDemand_preferIl_q01}
\eea

First, we will show that relations  (\ref{eqn_lm_MultiplicativeDemand_preferIl_q01}) are true at $t=T-1$.

For all $I\leq - q_T$ and $d_s\geq 0$, we have
\bean
J_T(I,d_s,0)&=&R_{s,T}(d_s) - \E H(I +q_T-\epsilon_{s,t}d_s) + \E [e\min\{I + q_T-\epsilon_{s,t}d_s,0\}]\\
&\leq & R_{s,T}(0) + (\e+\s)d_s  + \s(I + q_T-d_s) + \e(I + q_T-d_s)\\
&=& R_{s,T}(0) + \s(I+q_T) + \e(I+q_T)\\
&=& J_T(I,0, 0)
\eean
The inequality follows from $c_e > R_{s,T}^{\prime}(0) -\s \geq R_{s,T}^{\prime}(d_s) -\s$. Hence, $d_{s,T}^*(I)=0$ for all $I\leq - q_T$.

For every $I> - q_T$, we have that $\frac{\partial J_T}{\partial d_s}|_{d_s=0} = R_{s,T}^{\prime}(0) - \h >0$. Hence $d_{s,T}^*(I)>0$ for all $I> - q_T$. Therefore $I_{s,T}^*= -q_T$ and $V_T^{\prime}(0)\leq V_T^{\prime}(-q_T)=R_{s,T}^{\prime}(0)<R_{\ell,T}^{\prime}(0)+\s$.

Combining the above with the condition, $\frac{R_{\ell,T-1}^{\prime}(0)-\alpha R_{\ell,T}^{\prime}(0)}{\alpha} >\s$, we have $\frac{R_{\ell,T-1}^{\prime}(0)-\alpha R_{\ell,T}^{\prime}(0)}{\alpha}>c_p>V_T^{\prime}(0)-R_{\ell,T}^{\prime}(0)$, i.e., $\alpha V_T^{\prime}(0) < R_{\ell,T-1}^{\prime}(0)$. Hence, we can find a $\delta_2>0$ satisfying the following:
\bean
\alpha \frac{V_T(0) - V_T( -\overline{\overline{\epsilon_{s}}}\delta_2)}{\overline{\overline{\epsilon_s}}\delta_2}<\frac{R_{\ell,T-1}(\delta_2)-R_{\ell,T-1}(0)}{\delta_2}
\eean
From the concavity of $V_T$, we have that for any realization $\epsilon_{s,t}>0$
\bean
\alpha\frac{V_T(0) - V_T( - {\epsilon}_{s,t}\delta_2)}{ {\epsilon}_{s,t}\delta_2}
\leq \alpha \frac{V_T(0) - V_T( -\overline{\overline{\epsilon_{s}}}\delta_2)}{\overline{\overline{\epsilon_s}}\delta_2}
<\frac{R_{\ell,T-1}(\delta_2)-R_{\ell,T-1}(0)}{\delta_2}
\eean
\bean
\alpha\frac{V_T(0) - V_T( - {\epsilon}_{s,t}\delta_2)}{R_{\ell,T-1}(\delta_2)-R_{\ell,T-1}(0)}
<\frac{ {\epsilon}_{s,t}\delta_2}{\delta_2} \eean
If $ {\epsilon}_{s,T-1}=0$, then $\alpha\frac{V_T(0) - V_T( - {\epsilon}_{s,T-1}\delta_2)}{R_{\ell,T-1}(\delta_2)-R_{\ell,T-1}(0)}=\frac{ {\epsilon}_{s,T-1}\delta_2}{\delta_2}$. Taking the expectations of both sides of the above inequality, we have
\bean
\alpha\frac{V_T(0) - \E V_T( -\epsilon_{s,T-1}\delta_2)}{R_{\ell,T-1}(\delta_2)-R_{\ell,T-1}(0)}<\frac{\E(\epsilon_{s,T-1}\delta_2)}{\delta_2} =1.
\eean
\bea
 R_{\ell,T-1}(0) + \alpha V_T(0) < R_{\ell,T-1}(\delta_2) + \alpha \E V_T( -\epsilon_{s,T-1}\delta_2)
 \label{eqn_lm_MultiplicativeDemand_preferIl_q02}
\eea

Therefore,
\bean
J_{T-1}(-q_{T-1},0,0)&=&R_{\ell,T-1}(0)+R_{s,T-1}(0) -\E H(0) + \alpha \E V_{T}(0)\\
&<&R_{\ell,T-1}(\delta_2) + R_{s,T-1}(0) -\E H(0) + \alpha \E V_T( -\epsilon_{s,T-1}\delta_2) \\
&=& J_{T-1}(-q_{T-1},0,\delta_2).
\eean
The inequality follows from Ineq. (\ref{eqn_lm_MultiplicativeDemand_preferIl_q02}).

Hence, $d_{s,T-1}^*(-q_{T-1})+d_{\ell,T-1}^*(-q_{T-1})>0$. It is to say , $\min\{I_{\ell,t}^*,I_{s,t}^*\} + q_{T-1}< 0$.

Next we will show that  $I_{\ell,T-1}^*<I_{s,T-1}^*$ by contradiction. Suppose $I_{\ell,T-1}^*\geq I_{s,T-1}^*$, then $d_{s,T-1}^*(I_{s,T-1}^*) = d_{s,T-1}^*(I_{s,T-1}^*) =0$
\bean
\frac{\partial J_{T-1}}{\partial d_\ell}\big|_{I=I_{s,T-1}^*, d_s=0, d_\ell=0}
&=& R_{\ell,T-1}^{\prime}(0) - \alpha \E V_{T}^{\prime}(I)\leq 0.\\
\frac{\partial J_{T-1}}{\partial d_s}\big|_{I=I_{s,T-1}^*, d_s=0, d_\ell=0}
&=& R_{s,T-1}^{\prime}(0) + \E H^{\prime}(I + q_{T-1}) + \alpha \E V_{T}^{\prime}(I + q_{T-1})\\
&=& R_{s,T-1}^{\prime}(0) - \s + \alpha \E V_{T}^{\prime}(I + q_{T-1})=0.
\eean
Hence, $R_{s,T-1}^{\prime}(0)-\s\geq R_{\ell,T-1}^{\prime}(0)$ which  contradicts our assumption and therefore $I_{\ell,T-1}^*<I_{s,T-1}^*$. We still have that
\bean
V_{T-1}^{\prime}(0)\leq V_{T-1}^{\prime}(I_{\ell,T-1}^*)=R_{\ell,T-1}^{\prime}(0) - \E H^{\prime}(I_{s,T-1}^* + q_{T-1}) \leq R_{\ell,T-1}^{\prime}(0) +\s .
\eean
The first inequality follows from $I_{\ell,T-1}^*<-q_{T-1}\leq 0$ and the concavity of $V_{T-1}$.

Assume $V_{t+1}^{\prime}(0)\leq R_{\ell,t+1}^{\prime}(0) +\s$. We will show that $I_{\ell,t}^*<I_{s,t}^*$ and $V_t^{\prime}(0)\leq R_{\ell,t}^{\prime}(0) +s$. From $\frac{R_{\ell,t}^{\prime}(0)-\alpha R_{\ell,t+1}^{\prime}(0)}{\alpha} >\s$, we will have $\alpha V_{t+1}^{\prime}(0) < R_{\ell,t}^{\prime}(0)$. Hence, we can find a $\delta_3>0$ satisfying the inequality below
\bean
\alpha \frac{V_{t+1}(0) - V_{t+1}( -\overline{\overline{\epsilon_{s}}}\delta_3)}{\overline{\overline{\epsilon_s}}\delta_3}
<\frac{R_{\ell,t}(\delta_3)-R_{\ell,t}(0)}{\delta_3} .
\eean
Note    that for any realization with $ {\epsilon}_{s,t}=0$, we have:
\bea
 \alpha\frac{V_{t+1}(0) - V_{t+1}( - {\epsilon}_{s,t}\delta_3)}{R_{\ell,t}(\delta_3)-R_{\ell,t}(0)}=0=\frac{ {\epsilon}_{s,t}\delta_3}{\delta_3}
  \label{eqn-r0}
\eea
In addition,  for any realization $ {\epsilon}_{s,t}>0$ the following inequality holds.
\bea
\alpha\frac{V_{t+1}(0) - V_{t+1}( - {\epsilon}_{s,t}\delta_3)}{R_{\ell,t}(\delta_3)-R_{\ell,t}(0)}
<\frac{ {\epsilon}_{s,t}\delta_3}{\delta_3}.
 \label{eqn-r1}
\eea
Inequality (\ref{eqn-r1}), is due to the  the concavity of $V_{t+1}$, which implies:
\bean
\alpha\frac{V_{t+1}(0) - V_{t+1}( - {\epsilon}_{s,t}\delta_3)}{ {\epsilon}_{s,t}\delta_3}
\leq \alpha \frac{V_{t+1}(0) - V_{t+1}( -\overline{\overline{\epsilon_{s}}}\delta_3)}{\overline{\overline{\epsilon_s}}\delta_3}
<\frac{R_{\ell,t}(\delta_3)-R_{\ell,t}(0)}{\delta_3} .
\eean

Taking   expectations  inequalities (\ref{eqn-r0}),  (\ref{eqn-r1}) imply:
\bean
\alpha\frac{V_{t+1}(0) - \E V_{t+1}( -\epsilon_{s,t}\delta_3)}{R_{\ell,t}(\delta_3)-R_{\ell,t}(0)}<\frac{\E(\epsilon_{s,t}\delta_3)}{\delta_3} =1,
\eean
and
\bea
 R_{\ell,t}(0) + \alpha V_{t+1}(0) < R_{\ell,t}(\delta_3) + \alpha \E V_{t+1}( -\epsilon_{s,t}\delta_3) .
 \label{eqn_lm_MultiplicativeDemand_preferIl_q03}
\eea

Therefore,
\bean
J_{t}(-q_{t},0,0)&=&R_{\ell,t}(0)+R_{s,t}(0) -\E H(0) + \alpha \E V_{t+1}(0)\\
&<&R_{\ell,t}(\delta_3) + R_{s,t}(0) -\E H(0) + \alpha \E V_{t+1}( -\epsilon_{s,t}\delta_3) \\
&=& J_{t}(-q_{t},0,\delta_3).
\eean
The inequality follows from Ineq.  (\ref{eqn_lm_MultiplicativeDemand_preferIl_q03}).

Hence, $d_{s,t}^*(-q_{t})+d_{\ell,t}^*(-q_{t})>0$, i.e.,
$$\min\{I_{\ell,t}^*,I_{s,t}^*\} + q_{t}< 0 .$$

We next show that  $I_{\ell,t}^*<I_{s,t}^*$ by contradiction. Suppose $I_{\ell,t}^*\geq I_{s,t}^*$, Then $d_{s,t}^*(I_{s,t}^*) = d_{s,t}^*(I_{s,t}^*) =0$
\bean
\frac{\partial J_{t}}{\partial d_\ell}\big|_{I=I_{s,t}^*, d_s=0, d_\ell=0}
&=& R_{\ell,t}^{\prime}(0) - \alpha \E V_{t+1}^{\prime}(I)\leq 0.\\
\frac{\partial J_{t}}{\partial d_s}\big|_{I=I_{s,t}^*, d_s=0, d_\ell=0}
&=& R_{s,t}^{\prime}(0) + \E H^{\prime}(I + q_{t}) + \alpha \E V_{t+1}^{\prime}(I + q_{t})\\
&=& R_{s,t}^{\prime}(0) - \s + \alpha \E V_{t+1}^{\prime}(I + q_{t})=0.
\eean
Hence, $R_s^{\prime}(0)-\s\geq R_\ell^{\prime}(0)$ which  contradicts our assumption. Hence, we have $I_{\ell,t}^*<I_{s,t}^*$. We still have that
\bean
V_{t}^{\prime}(0)\leq V_{t}^{\prime}(I_{\ell,t}^*)=R_{\ell,t}^{\prime}(0) - \E H^{\prime}(I_{s,t}^* + q_{t}) \leq R_{\ell,t}^{\prime}(0) +\s .
\eean
The first inequality follows from $I_{\ell,t}^*<-q_{t}\leq 0$ and the concavity of $V_{t}$. Hence, we obtain (\ref{eqn_lm_MultiplicativeDemand_preferIl_q01}). The proof is now complete.
\end{pf}

\begin{pf}{\bf of Theorem~\ref{lem_MultiplicativeDemand_samefuntion_IsIl}.} We need to consider the following three cases:

\emph{a)} $d^*_{s,t}(I)\leq d^*_{\ell,t}(I)$ when  $ \E [ \epsilon_{s,t}H^{\prime}(I + q_t - \epsilon_{s,t} d^*_{s,t}(I))] <0$.

\emph{b)} $d^*_{s,t}(I)\geq d^*_{\ell,t}(I)$ when  $ \E [ \epsilon_{s,t}H^{\prime}(I + q_t - \epsilon_{s,t} d^*_{s,t}(I))] >0$.

\emph{c)} $d^*_{s,t}(I) = d^*_{\ell,t}(I)$ when  $ \E [ \epsilon_{s,t}H^{\prime}(I + q_t - \epsilon_{s,t} d^*_{s,t}(I))] =0$.

a) We  prove case a) by contradiction. Suppose $d^*_{s,t}(I)> d^*_{\ell,t}(I)$. If $\E [ \epsilon_{s,t}H^{\prime}(I + q_t - \epsilon_{s,t} d^*_{s,t}(I))] < 0$, then we can find $\delta_1 \in (0, d^*_{s,t}(I) - d^*_{\ell,t}(I))$ satisfying that
\bea
\E [ H(I + q_t - \epsilon_{s,t} (d^*_{s,t}(I)-\delta_1))] < \E [ H(I + q_t - \epsilon_{s,t} d^*_{s,t}(I))].\label{eqn_lm_MultiplicativeDemand_samefuntion_IsIl_q1}
\eea
From the concavity of $R_{s,t}$ and $d^*_{\ell,t}(I) + \delta_1 < d^*_{s,t}(I)$, we have
\bean
R_{s,t}(d^*_{s,t}(I)) - R_{s,t}(d^*_{s,t}(I) - \delta_1) < R_{s,t}(d^*_{\ell,t}(I) + \delta_1) - R_{s,t}(d^*_{\ell,t}(I)) =  R_{\ell,t}(d^*_{\ell,t}(I) + \delta_1) - R_{\ell,t}(d^*_{\ell,t}(I)) .
\eean
That is to say,
\bea
R_{s,t}(d^*_{s,t}(I)) + R_{\ell,t}(d^*_{\ell,t}(I))  < R_{s,t}(d^*_{s,t}(I) - \delta_1) + R_{\ell,t}(d^*_{\ell,t}(I) + \delta_1) .
\label{eqn_lm_MultiplicativeDemand_samefuntion_IsIl_q2}
\eea
Meanwhile, from the concavity of $V_{t+1}$ and $\delta_1 < d^*_{s,t}(I) - d^*_{\ell,t}(I)$, we have
\bea
&&\E V_{t+1}(I + q_t - \epsilon_{s,t} d^*_{s,t}(I) - \epsilon_{\ell,t} d^*_{\ell,t}(I))\nonumber\\
&=& (1- \frac{\delta_1}{d^*_{s,t}(I) - d^*_{\ell,t}(I)}) \E V_{t+1}(I + q_t - \epsilon_{s,t} d^*_{s,t}(I) - \epsilon_{\ell,t} d^*_{\ell,t}(I)) \nonumber\\
&& \,\,\,\, +
\frac{\delta_1}{d^*_{s,t}(I) - d^*_{\ell,t}(I)} \E V_{t+1}(I + q_t - \epsilon_{s,t} d^*_{\ell,t}(I) - \epsilon_{\ell,t} d^*_{s,t}(I))\nonumber\\
&\leq & \E V_{t+1}(I + q_t - \epsilon_{s,t} (d^*_{s,t}(I)-\delta_1) - \epsilon_{\ell,t} (d^*_{\ell,t}(I)+\delta_1)) \label{eqn_lm_MultiplicativeDemand_samefuntion_IsIl_q3}
\eea
Together with Ineqs.  (\ref{eqn_lm_MultiplicativeDemand_samefuntion_IsIl_q1}) and (\ref{eqn_lm_MultiplicativeDemand_samefuntion_IsIl_q2}), we have
\bean
&&J_t(I,d^*_{s,t}(I)-\delta_1, d^*_{\ell,t}(I)+\delta_1)\\
&=& R_{s,t}(d^*_{s,t}(I)-\delta_1) + R_{\ell,t}(d^*_{\ell,t}(I)+\delta_1) - \E [ H(I + q_t - \epsilon_{s,t} (d^*_{s,t}(I)-\delta_1))] \\
&& \,\,\,\,\, + \alpha \E V_{t+1}(I + q_t - \epsilon_{s,t} (d^*_{s,t}(I)-\delta_1) - \epsilon_{\ell,t} (d^*_{\ell,t}(I)+\delta_1)) \\
&>& R_{s,t}(d^*_{s,t}(I)) + R_{\ell,t}(d^*_{\ell,t}(I)) - \E [ H(I + q_t - \epsilon_{s,t} d^*_{s,t}(I))] \\
&& \,\,\,\,\, + \alpha \E V_{t+1}(I + q_t - \epsilon_{s,t} d^*_{s,t}(I) - \epsilon_{\ell,t} d^*_{\ell,t}(I))\\
&=&J_t(I,d^*_{s,t}(I), d^*_{\ell,t}(I))
\eean
which contradicts with the optimality of $(d_{s,t}^*(I),d_{\ell,t}^*(I))$.

b) The  proof  of case b) is also by contradiction. Suppose $d^*_{s,t}(I)< d^*_{\ell,t}(I)$. If $\E [ \epsilon_{s,t}H^{\prime}(I + q_t - \epsilon_{s,t} d^*_{s,t}(I))] > 0$, then we can find $\delta_2 \in (0, d^*_{\ell,t}(I) - d^*_{s,t}(I))$ satisfying that
\bea
\E H(I + q_t - \epsilon_{s,t} (d^*_{s,t}(I)+\delta_2)) < \E H(I + q_t - \epsilon_{s,t} d^*_{s,t}(I))\label{eqn_lm_MultiplicativeDemand_samefuntion_IsIl_q4}
\eea
From the concavity of $R_\ell$ and $d^*_{s,t}(I) +  \delta_2 < d^*_{\ell,t}(I)$, we have
\bean
R_{\ell,t}(d^*_{\ell,t}(I)) - R_{\ell,t}(d^*_{\ell,t}(I) - \delta_2) < R_{\ell,t}(d^*_{s,t}(I) + \delta_2) - R_{\ell,t}(d^*_{s,t}(I)) =  R_{s,t}(d^*_{s,t}(I) + \delta_2) - R_{s,t}(d^*_{s,t}(I))
\eean
That is to say,
\bea
R_{s,t}(d^*_{s,t}(I)) + R_{\ell,t}(d^*_{\ell,t}(I))  < R_{\ell,t}(d^*_{\ell,t}(I) - \delta_2) + R_{s,t}(d^*_{s,t}(I) + \delta_2)\label{eqn_lm_MultiplicativeDemand_samefuntion_IsIl_q5}
\eea
Meanwhile, from the concavity of $V_{t+1}$ and $\delta_2 < d^*_{\ell,t}(I) - d^*_{s,t}(I)$, we have:
\bea
&&\E V_{t+1}(I + q_t - \epsilon_{s,t} d^*_{s,t}(I) - \epsilon_{\ell,t} d^*_{\ell,t}(I))\nonumber\\
&=& (1- \frac{\delta_2}{d^*_{\ell,t}(I) - d^*_{s,t}(I)}) \E V_{t+1}(I + q_t - \epsilon_{s,t} d^*_{s,t}(I) - \epsilon_{\ell,t} d^*_{\ell,t}(I)) \nonumber\\
&& \,\,\,\, +
\frac{\delta_2}{d^*_{\ell,t}(I) - d^*_{s,t}(I)} \E V_{t+1}(I + q_t - \epsilon_{s,t} d^*_{\ell,t}(I) - \epsilon_{\ell,t} d^*_{s,t}(I))\nonumber\\
&\leq & \E V_{t+1}(I + q_t - \epsilon_{s,t} (d^*_{s,t}(I)+\delta_2) - \epsilon_{\ell,t} (d^*_{\ell,t}(I)-\delta_2)) \label{eqn_lm_MultiplicativeDemand_samefuntion_IsIl_q6}
\eea
Together with (\ref{eqn_lm_MultiplicativeDemand_samefuntion_IsIl_q4}) and (\ref{eqn_lm_MultiplicativeDemand_samefuntion_IsIl_q5}), we have
\bean
&&J_t(I, d^*_{s,t}(I)+\delta_2, d^*_{\ell,t}(I)-\delta_2)\\
&=& R_{s,t}(d^*_{s,t}(I)+\delta_2) + R_{\ell,t}(d^*_{\ell,t}(I)-\delta_2) - \E [ H(I + q_t - \epsilon_{s,t} (d^*_{s,t}(I)+\delta_2))] \\
&& \,\,\,\,\, + \alpha \E V_{t+1}(I + q_t - \epsilon_{s,t} (d^*_{s,t}(I)+\delta_2) - \epsilon_{\ell,t} (d^*_{\ell,t}(I)-\delta_2)) \\
&>& R_{s,t}(d^*_{s,t}(I)) + R_{\ell,t}(d^*_{\ell,t}(I)) - \E [ H(I + q_t - \epsilon_{s,t} d^*_{s,t}(I))] \\
&& \,\,\,\,\, + \alpha \E V_{t+1}(I + q_t - \epsilon_{s,t} d^*_{s,t}(I) - \epsilon_{\ell,t} d^*_{\ell,t}(I))\\
&=&J_t(I,d^*_{s,t}(I), d^*_{\ell,t}(I))
\eean
and this also leads to a contradiction with the optimality of $(d_{s,t}^*(I),d_{\ell,t}^*(I))$.

c) To prove case c), we simply point to the continuity of $d_{s,t}^*(I), d_{\ell,t}^*(I)$ and the results of cases \emph{a)} and \emph{b)}.

\end{pf}

\begin{pf} {\bf of Proposition~\ref{lem_MultiplicativeDemand_Correlated_positivedst}.}
When $d_{s,t}^*(I)>0$, we have for every realization $\tilde \epsilon_{s,t}$
\bea
V_t(I) -V_t(I-\tilde \epsilon_{s,t}d_{s,t}^*(I))
\leq V_t(I) - J_t(I-\tilde \epsilon_{s,t} d_{s,t}^*(I),0,d_{\ell,t}^*(I))\label{eqn_lem_MultiplicativeDemand_Correlated_positivedst_q1}
\eea
Hence, we have
\bean
V_t(I) -V_t(I-d_{s,t}^*(I))&\leq& V_t(I) - \E_{\tilde\epsilon_{s,t}}[V_t(I-\tilde \epsilon_{s,t}d_{s,t}^*(I))]\\
&\leq& V_t(I) - \E_{\tilde\epsilon_{s,t}}[J_t(I-\tilde \epsilon_{s,t}d_{s,t}^*(I),0,d_{\ell,t}^*(I))]\\
&=& R_{s,t}(d_{s,t}^*(I)) + R_{\ell,t}(d_{\ell,t}^*(I)) - \E H(I + q_t -\epsilon_{s,t}d_{s,t}^*(I))\\
&& \, \, \, \, + \alpha \E V_{t+1}(I + q_t -\epsilon_{s,t}d_{s,t}^*(I)-\epsilon_{\ell,t}d_{\ell,t}^*(I))  - \E_{\tilde\epsilon_{s,t}}J_t(I-\tilde \epsilon_{s,t}d_{s,t}^*(I),0,d_{\ell,t}^*(I))\\
&=&  R_{s,t}(d_{s,t}^*(I)) -R_{s,t}(0)
\eean
The second inequality follows from Jensen's inequality. The second inequality follows from (~\ref{eqn_lem_MultiplicativeDemand_Correlated_positivedst_q1}).
\bean
V_t^{\prime}(I)\leq \frac{V_t(I) -V_t(I-d_{s,t}^*(I))}{d_{s,t}^*(I)}
\leq \frac{R_{s,t}(d_{s,t}^*(I)) - R_{s,t}(0)}{d_{s,t}^*(I)}
<  R_{s,t}^{\prime}(0).
\eean
\end{pf}

\begin{pf} {\bf of Proposition~\ref{lem_MultiplicativeDemand_perfectcorrelated}.} If $I+q_t\leq 0$, then from the first order partial derivatives  we have:
Notice that
\bean
0&=&\frac{\partial J_t(I,d_s,d_\ell)}{\partial d_s}\bigg|_{d_s=0} = R_{s,t}^{\prime}(0)  + \E\big[ \epsilon_{s,t} (-c_p)\big] - \alpha \E \big[\epsilon_{s,t} V_{t+1}^{\prime}(I + q_t  - \epsilon_{\ell,t} d_{\ell,t}^*(I))\big]\\
&=& R_{s,t}^{\prime}(0)  - \big\{ c_p +\alpha \E \big[\epsilon_{s,t} V_{t+1}^{\prime}(I + q_t  - \epsilon_{\ell,t} d_{\ell,t}^*(I))\big]\}.
\eean
\bean
V^{\prime}(I)=\frac{\partial J_t(I,d_s,d_\ell)}{\partial I}=\s + \E V_{t+1}^{\prime}(I + q_t  - \epsilon_{\ell,t} d_{\ell,t}^*(I)).
\eean
Notice that if $\epsilon_{s,t}$ and $\epsilon_{\ell,t}$ are perfect positive correlated, then from the concavity of $V(.)$ we have
\bean
\E V_{t+1}^{\prime}(I + q_t  - \epsilon_{\ell,t} d_{\ell,t}^*(I))=\E \big[\epsilon_{s,t}\big]\E V_{t+1}^{\prime}(I + q_t  - \epsilon_{\ell,t} d_{\ell,t}^*(I))
\leq \E \big[\epsilon_{s,t} V_{t+1}^{\prime}(I + q_t  - \epsilon_{\ell,t} d_{\ell,t}^*(I))\big].
\eean
If $\epsilon_{s,t}$ and $\epsilon_{\ell,t}$ are perfect negative correlated, then from the concavity of $V(.)$ we have
\bean
\E V_{t+1}^{\prime}(I + q_t  - \epsilon_{\ell,t} d_{\ell,t}^*(I))= \E \big[\epsilon_{s,t}\big] \E V_{t+1}^{\prime}(I + q_t  - \epsilon_{\ell,t} d_{\ell,t}^*(I))
\geq \E \big[\epsilon_{s,t} V_{t+1}^{\prime}(I + q_t  - \epsilon_{\ell,t} d_{\ell,t}^*(I))\big].
\eean
Therefore $R_{s,t}^{\prime}(0)-V^{\prime}(I)\geq 0$ if  $\epsilon_{s,t}$ and $\epsilon_{\ell,t}$ perfect positive correlated and  $R_{s,t}^{\prime}(0)-V^{\prime}(I)\leq 0$ if  $\epsilon_{s,t}$ and $\epsilon_{\ell,t}$ perfect negative correlated

Follow the similar logic, we can obtain the statement is true when $I+q_t>0$.  We conclude the  proof.
\end{pf}

\begin{pf}{\bf of Theorem~\ref{th_perfectnegativecorrelated_Chooseingon-siteStore_Threshold}.}
The result follows from Proposition~\ref{lem_MultiplicativeDemand_Correlated_positivedst} and Proposition~\ref{lem_MultiplicativeDemand_perfectcorrelated}.
\end{pf}

\begin{pf}{\bf of Lemma~\ref{lem_SS_Unique_Price}.}
From the concavity of $V_t$ and the linear combination form of the  demand functions, we have $J_t^u$ is submodular in $(I,p)$. Hence, $p^u_t(I)$ decreases in $I$
\end{pf}

\begin{pf}{\bf of Theorem~\ref{th_SS_Unique_Sandwich}.} To simplify the proof, we denote
\bean
G(I)&=&-\beta_{s,t}[R_{s,t}^{\prime}(0) + H^{\prime}(I+q_t) - V_{t+1}^{\prime}(I+q_t)]\\
&& \,\,\,\,\, \,\,\, -\beta_{\ell,t}[R_{\ell,t}^{\prime}(0) - V_{t+1}^{\prime}(I+q_t)].
\eean
Notice that $G(I)$ decreases in $I$.

Taking first order partial derivatives of $J_t^u (I, p)$ with respect of $p$,
\bean
\frac{dJ_t^u}{dp}\bigg|_{I=I_t^u*}&=& -\beta_{s,t}R_{s,t}^{\prime}(0)  - \beta_{\ell,t}R_{\ell,t}^{\prime}(0) - \E\big[\varepsilon_{s,t}\beta_{s,t}H^{\prime}(I+q_t)\big]\\
&&\,\,\, + \E\big[(\varepsilon_{s,t}\beta_{s,t}+\varepsilon_{\ell,t}\beta_{\ell,t})V_{t+1}^{\prime}(I+q_t)\big]\\
&=& G(I_t^u*)
\eean

To compare $I_{s,t}^*$, $I_{\ell,t}^*$ and $I_t^u*$, we discuss two cases: (i) $I_{s,t}^*\geq I_{f,t}^*$; (ii) $I_{s,t}^*\leq I_{f,t}^*$.

To see (i), we need to prove that $I_{\ell,t}^*\leq I_t^u* \leq I_{s,t}^*$. According to first order partial derivatives of $J_t$ with respect to $d_s$ and $d_l$.
\bean
&0=\frac{\partial J_t}{d_s}\bigg|_{I=I_{f,t}^*}
\leq R_{s,t}^{\prime}(0) + H^{\prime}(I+q_t) - V_{t+1}^{\prime}(I+q_t), \,\,\, \,
0=\frac{\partial J_t}{d_l}\bigg|_{I=I_{f,t}^*}
\leq R_{\ell,t}^{\prime}(0) - V_{t+1}^{\prime}(I+q_t)\\
&0=\frac{\partial J_t}{d_s}\bigg|_{I=I_{s,t}^*}
= R_{s,t}^{\prime}(0) + H^{\prime}(I+q_t) - V_{t+1}^{\prime}(I+q_t), \,\,\, \,
0 \geq \frac{\partial J_t}{d_l}\bigg|_{I=I_{s,t}^*}
= R_{\ell,t}^{\prime}(0) - V_{t+1}^{\prime}(I+q_t).
\eean
Hence, $G(I_{\ell,t}^*)\geq 0$ and $G(I_{s,t}^*)\leq 0$. From the decreasing proposition of $G$, we have $I_{\ell,t}^*\leq I_t^u* \leq I_{s,t}^*$.

Following a similar argument, this result is also true for Case (ii) and we conclude the proof.
\end{pf}


\begin{thebibliography}{99}
\bibitem{Agatz2008}
Agatz, N. A., Fleischmann, M.,  Van Nunen, J. A. 2008. 
\newblock E-fulfillment and multi-channel distribution - A review. 
\newblock {\em European Journal of Operational Research,} 187(2), 339-356.

\bibitem{AgrawalSeshadri2000}
Agrawal, V.,  and Seshadri, S. 2000.
\newblock Impact of uncertainty and risk aversion on price and order quantity
  in the newsvendor problem.
\newblock {\em Manufacturing and Service Operations Management}, 2(4):410--422.

\bibitem{annreport}
{Ann Taylor}.
\newblock Ann taylor, 2012.  Annual report.
\newblock
  \url{http://investor.anninc.com/phoenix.zhtml?c=78167&p=irol-reportsannual},
  2012.

\bibitem{AnupindiAkella1993}
Anupindi, R., and R.~Akella. 1993.
\newblock Diversification under supply uncertainty.
\newblock {\em Management Science}, 39(8):944--963.

\bibitem{araman2009dynamic}
 Araman V.F.,  and R. Caldentey. 2009.
\newblock Dynamic pricing for nonperishable products with demand learning.
\newblock {\em Operations research}, 57(5):1169--1188.

\bibitem{arslan2007single}
Arslan, H., Graves, S.C., and T.~Roemer. 2007.
\newblock A single-product inventory model for multiple demand classes.
\newblock {\em Management Science}, 53(9):1486--1500.

\bibitem{AvivFedergruen2001a}
 Aviv, Y., and A. Federgruen. 2001.
\newblock Capacitated multi-item inventory systems with random and seasonally
  fluctuating demands: Implications for postponement strategies.
\newblock {\em Management Science}, 47(4):512--531.

\bibitem{AvivFedergruen2001b}
 Aviv, Y., and A. Federgruen. 2001.
\newblock Design for postponement: A comprehensive characterization of its
  benefits under unknown demand distributions.
\newblock {\em Operations Research}, 49(4):578--598.

\bibitem{barankin1961delivery}
Barankin, E.W. 1961.
\newblock A delivery-lag inventory model with an emergency provision (the
  single-period case).
\newblock {\em Naval Research Logistics Quarterly}, 8(3):285--311.

\bibitem{wsj2014}
 Bensinger, G. 2014.
\newblock Spreading Black Friday fever to China's shoppers.
\newblock
  \url{http://www.wsj.com/articles/spreading-black-friday-fever-to-chinas-shoppers-1416954943},
  2014.
\newblock Accessed: 2014-12-11, Nov. 25.

\bibitem{BernsteinFedergruen2004}
Bernstein, F.,  and   Federgruen, A. 2004.
\newblock A general equilibrium model for industries with price and service
  competition.
\newblock {\em Operations Research}, 52(6):868--886.

\bibitem{BernsteinFedergruen2007}
Bernstein, F.,  and   Federgruen, A.  2007.
\newblock Coordination mechanisms for supply chains under price and service
  competition.
\newblock {\em Manufacturing and Service Operations Management}, 9(3):242--262.

\bibitem{bitran2003commissioned}
Bitran, G., and R.~Caldentey, R. 2003.
\newblock Commissioned paper: An overview of pricing models for revenue
  management.
\newblock {\em Manufacturing \& Service Operations Management}, 5(3):203--229.

\bibitem{bitran1998coordinating}
Bitran, G., Caldentey,  R., and S.~Mondschein. 1998.
\newblock Coordinating clearance markdown sales of seasonal products in retail
  chains.
\newblock {\em Operations research}, 46(5):609--624.

\bibitem{retailnet2012}
 Bomberowitz, J. 2012. 
\newblock Ship from store: What retailers can learn from the early adopters.
\newblock
  \url{http://connected.retailnetgroup.com/index.php/2012/09/24/ship-from-store-what-retailers-can-learn-from-the-early-adopters/}.

\bibitem{BrynjolfssonHuRahman2009}
 Brynjolfsson, E.,  Hu, Y., and M.~S. Rahman. 2009.
\newblock Battle of the retail channels: How product selection and geography
  drive cross-channel competition.
\newblock {\em Management Science}, 55(11):1755--1765.

\bibitem{BrynjolfssonHuSmith2003}
Brynjolfsson, E.,  Hu, Y., and M.~D. Smith. 2003.
\newblock Consumer surplus in the digital economy: Estimating the value of
  increased product variety at online booksellers.
\newblock {\em Management Science}, 49(11):1580--1596.

\bibitem{BrynjolfssonSmith2000}
Brynjolfsson E., and M.~D. Smith. 2009.
\newblock Battle of the retail channels: How product selection and geography
  drive cross-channel competition.
\newblock {\em Management Science}, 46(4):563--585.

\bibitem{Butler2005}
  Butler, E. 2005.
\newblock Orchestrating sept. 11 comeback: J\&r music finally regains lost
  sales with online focus, new product mix.
\newblock {\em Crain's New York Business}, Jan. 17.

\bibitem{CampbellFrei2010}
Campbell, D., and Frei, F. 2010.
\newblock Cost structure, customer profitability, and retention implications of
  self-service distribution channels: Evidence from customer behavior in an
  online banking channel.
\newblock {\em Management Science}, 56(1):4--24.

\bibitem{Chen2001}
Chen, F. 2001. 
\newblock Market segmentation, advanced demand information, and supply chain
  performance.
\newblock {\em Manufacturing and Service Operations Management}, 3(1):53--67.
  

\bibitem{Chen21092008}
Chen, K.Y., Kaya, M., and O. Ozer. 2008.
\newblock Dual sales channel management with service competition.
\newblock {\em Manufacturing and Service Operations Management},
  10(4):654--675.

\bibitem{Chen2013}
Chen, W. Feng, Q., and S.~Seshadri. 2013.
\newblock Sourcing from suppliers with random yield for price-dependent demand.
\newblock {\em Annals of Operations Research}, 208(1):557--579.

\bibitem{ChenSimchi-Levi2004a}
 Chen, X., and   Simchi-Levi, D. 2004.
\newblock Coordinating inventory control and pricing strategies with random
  demand and fixed ordering cost: The finite horizon case.
\newblock {\em Operations Research}, 52(6):887--896.

\bibitem{ChenSimchi-Levi2004b}
 Chen, X., and    Simchi-Levi, D. 2004.
\newblock Coordinating inventory control and pricing strategies with random
  demand and fixed ordering cost: the infinite horizon case.
\newblock {\em Mathematics of Operations Research}, 29(3):698--723.

 
\bibitem{ChenSimchi-Levi2006}
Chen, X., and   Simchi-Levi, D. 2006.
\newblock Coordinating inventory control and pricing strategies with random
  demand and fixed ordering cost: the continuous review model.
\newblock {\em Operations Research Letters}, 34:323--332.

\bibitem{cohen1977joint}
Cohen, M.A. 1977.
\newblock Joint pricing and ordering policy for exponentially decaying
  inventory with known demand.
\newblock {\em Naval Research Logistics Quarterly}, 24(2):257--268.

\bibitem{datta2001}
Datta, T.K., and P. Karabi. 2001
\newblock An inventory system with stock-dependent, price-sensitive demand
  rate.
\newblock {\em Production planning \& control}, 12(1):13--20.

\bibitem{dave1982probabilistic}
Dave, U., and Y.K.~Shah. 1982.
\newblock A probabilistic inventory model for deteriorating items with lead
  time equal to one scheduling period.
\newblock {\em European Journal of Operational Research}, 9(3):281--285.

\bibitem{dong2009dynamic}
  Dong, L.,  Kouvelis, P., and Z. Tian. 2009.
\newblock Dynamic pricing and inventory control of substitute products.
\newblock {\em Manufacturing \& Service Operations Management}, 11(2):317--339.

\bibitem{elmaghraby2003dynamic}
Elmaghraby, W., and P.~Keskinocak. 2003.
\newblock Dynamic pricing in the presence of inventory considerations: Research
  overview, current practices, and future directions.
\newblock {\em Management Science}, pages 1287--1309.

\bibitem{FedergruenHeching1999}
Federgruen, A., and A.~Heching. 1999.
\newblock Combined pricing and inventory control under uncertainty.
\newblock {\em Operations Research}, pages 454--475.

\bibitem{federgruen2002multilocation}
Federgruen, A., and A.~Heching. 2002.
\newblock Multilocation combined pricing and inventory control.
\newblock {\em Manufacturing \& Service Operations Management}, 4(4):275--295.

\bibitem{FedYang2009}
Federgruen, A., and   N. Yang. 2011. 
\newblock Procurement strategies with unreliable suppliers.
\newblock {\em Operations Research}, 59(4):1033--1039.

\bibitem{feng2000continuous}
Feng, Y., and B.~Xiao. 2000.
\newblock A continuous-time yield management model with multiple prices and
  reversible price changes.
\newblock {\em Management Science}, 46(5): 644--657.

\bibitem{FormanGhoseGoldfarb2009}
  Forman, C. Ghose, A., and A. Goldfarb.  2009.
\newblock Competition between local and electronic markets: How the benefit of
  buying online depends on where you live.
\newblock {\em Management Science}, 55:47--57.

\bibitem{fukuda1964optimal}
Fukuda, Y. 1964.
\newblock Optimal policies for the inventory problem with negotiable leadtime.
\newblock {\em Management Science}, pages 690--708.

\bibitem{GallegoRyzin1994}
Gallego, G., and G.~J. van Ryzin. 1994.
\newblock Optimal dynamic pricing of inventories with stochastic demand over
  finite horizons.
\newblock {\em Management Science}, 40(8):999--1020.

\bibitem{Heching2002}
Heching, A.,  Gallego, G., and  van Ryzin G. 2002. 
\newblock Mark-down pricing: An empirical analysis of policies and revenue potential at one apparel retailer.
\newblock {\em Journal of Revenue and Pricing Management},  1(2): 139--160.

\bibitem{HeymanSobel1984}
Heyman, D.~P.,  and M.~J. Sobel. 1984.
\newblock {\em Stochastic Models in Operations Research, Vol. II.}
\newblock McGraw-Hill, New York.

\bibitem{Huang2009}
Huang, W., and J.~M. Swaminathan. 2009.
\newblock Introduction of a second channel: Implications for pricing and
  profits.
\newblock {\em European Journal of Operational Research}, 194(1):258--279.

\bibitem{PentinaPeltonHasty2009}
Pelton. L.~E., Pentina, I., and R.~W. Hasty. 2009.
\newblock Performance implications of online entry timing by store-based
  retailers: A longitudinal investigation.
\newblock {\em Journal of Retailing}, 85(2):177--193.

\bibitem{karlin1958inventory}
Karlin, S., and H.~Scarf. 1958.
\newblock Inventory models of the arrow-harris-marschak type with time lag.
\newblock {\em Studies in the mathematical theory of inventory and production},
  pages 155--178.

\bibitem{qr2012}
Katehakis, M.N. and  L.C. Smit. 2012.
\newblock On computing optimal (q, r) replenishment policies under quantity
  discounts.
\newblock {\em Annals of Operations Research} {200}(1) 279--298.


\bibitem{Kazaz2004}
Kazaz, B.  2004.
\newblock Production planning under yield and demand uncertainty with
  yield-dependent cost and price.
\newblock {\em Manufacturing and Service Operations Management}, 6(3):209--224.
 

\bibitem{Khermouch2012}
 Khermouch, G. 2012.
\newblock Retailer of the year.
\newblock {\em Beverage Dynamics}.
\newblock January/February, available at
  \texttt{www.garyswine.com/retailer\_year.pdf}.

\bibitem{levin2008risk}
Levin Y., McGill, J., and M.~Nediak. 2008.
\newblock Risk in revenue management and dynamic pricing.
\newblock {\em Operations Research}, 56(2):326.

\bibitem{Li01021992}
 Li, L. 1992.
\newblock The role of inventory in delivery-time competition.
\newblock {\em Management Science}, 38(2):182--197.

\bibitem{LiZheng2006}
Li, Q., and S.~Zheng. 2006.
\newblock Joint inventory replenishment and pricing control for systems with
  uncertain yield and demand.
\newblock {\em Operations Research}, 54:696--705.

\bibitem{MonahanPetruzziZhao2004}
Monahan, G.E,  Petruzzi, N.C., and W.~Zhao. 2004.
\newblock The dynamic pricing problem from a newsvendor's perspective.
\newblock {\em Manufacturing and Service Operations Management}, 6(1):73--91.

\bibitem{pang2012technical}
 Pang, Z.  Chen, F.Y., and Y. Feng. 2012.
\newblock Technical note—a note on the structure of joint inventory-pricing
  control with leadtimes.
\newblock {\em Operations Research}, 60(3):581--587.

\bibitem{Petruzzi1999}
Petruzzi, N.C., and M. Dada. 1999.
\newblock Pricing and the newsvendor problem: A review with extensions.
\newblock {\em Operations Research}, 47(2):183--194.

\bibitem{martin2013}
Shi, J., Katehakis,  M.N., and B.~Melamed. 2013.
\newblock Martingale methods for pricing inventory penalties under continuous
  replenishment and compound renewal demands.
\newblock In M.N. Katehakis, S.M. Ross, and J.~Yang, editors, {\em Cyrus Derman
  Memorial Volume I: Optimization under Uncertainty: Costs, Risks and
  Revenues}, volume 208, pages 593--612. Annals of Operations Research,
  Springer, New York.

\bibitem{pis2014}
Shi, J., Katehakis,  M.N.,  Melamed, B., and Y. Xia. 2014. 
\newblock Production-inventory systems with lost sales and compound poisson
  demands.
\newblock {\em Operations Research}, 6(5):1048 -- 1063.

\bibitem{simchi2005logic}
Simchi-Levi, D.,  Bramel, J., and X.~Chen. 2005.
\newblock {\em The logic of logistics: theory, algorithms, and applications for
  logistics and supply chain management}.
\newblock Springer Verlag.

\bibitem{SmithAchabal1998}
Smith, S.~A.,  and D.~D. Achabal. 1998.
\newblock Clearance pricing and inventory policies for retail chains.
\newblock {\em Management Science}, 44(3):285--300.

\bibitem{song2009technical}
 Song, Y., Ray, S., and T. Boyaci. 2009.
\newblock Technical note—optimal dynamic joint inventory-pricing control for
  multiplicative demand with fixed order costs and lost sales.
\newblock {\em Operations Research}, 57(1):245--250.
\bibitem{bl}
Spencer S. 2015.
\newblock China Shops Alibaba for U.S. Goods From Toothbrushes to Nuts
\newblock {\em Bloomberg Business} 
April 21, 2015 
\url{http://www.bloomberg.com/news/articles/2015-04-22/china-shops-alibaba-for-u-s-goods-from-toothbrushes-to-nuts}
\bibitem{wsj2013}
 Stevens, L. 2013.
\newblock Retailers turn store clerks into web shippers.
\newblock
  \url{http://www.wsj.com/articles/SB10001424052702303332904579228602337333952},
  Dec. 9 . 

\bibitem{talluri2005theory}
Talluri, K.T.,  and G.~Van~Ryzin. 2005.
\newblock {\em The theory and practice of revenue management}, volume~68.
\newblock Springer Verlag.

\bibitem{Thomas1974}
 Thomas, L.J. 1974.
 \newblock Technical Note?Price and Production Decisions with Random Demand. \newblock {\em Operations Research} 22(3):513-518.

\bibitem{intRetail2014}
 Tong F. 2014.
\newblock China's big online marketplaces will ship to overseas {C}hinese.
\newblock
 www.internetretailer.com/2013/07/26/chinas-big-online-marketplaces-will-ship-overseas-chinese.


\bibitem{birbil2010tractable}
 Topaloglu, H.,  Birbil, S.I.  Frenk, J.B.G., and N.  Noyan. 2012.
\newblock Tractable open loop policies for joint overbooking and capacity
  control over a single flight leg with multiple fare classes.
\newblock {\em Transportation Science}, 46(4):460--481.

\bibitem{Topkis1998}
Topkis, D.~M.  1998.
\newblock {\em Supermodularity and Complementarity}.
\newblock Princeton University Press, Princeton NJ.

\bibitem{veinott1966status}
Veinott A.F.~Jr . 1966.
\newblock The status of mathematical inventory theory.
\newblock {\em Management Science}, 12(1): 745--777.

\bibitem{wang2009inventory}
Wang, H., and H.~Yan. 2009.
\newblock Inventory management for customers with alternative lead times.
\newblock {\em Production and Operations Management}, 18(6):705--720.

\bibitem{wright1968optimal}
Wright, G.P. 1968.
\newblock Optimal policies for a multi-product inventory system with negotiable
  lead times.
\newblock {\em Naval Research Logistics Quarterly}, 15(3):375--401.

\bibitem{xu2005multi}
Xu, N..  2005.
\newblock Multi-period dynamic supply contracts with cancellation.
\newblock {\em Computers and {O}perations {R}esearch}, 32(12):3129--3142.

\bibitem{Yang2007}
 Yang, B.,  and   Geunes, J. 2007.
\newblock Inventory and lead time planning with lead-time-sensitive demand.
\newblock {\em IIE Transactions}, 39(5):439--452.

\bibitem{Yin2007113}
 Yin, R., and K. Rajaram.   2007.
\newblock Joint pricing and inventory control with a markovian demand model.
\newblock {\em European Journal of Operational Research}, 182(1):113 -- 126.

\bibitem{You1999}
 You, P.~S. 1999.
\newblock Dynamic pricing in airline seat management for flights with multiple
  legs.
\newblock {\em Transportation Science}, 34(2):192--206.

\bibitem{liu2006}
Yunchuan, L., Gupta, S., and Z.~J. Zhang. 2006.
\newblock Note on self-restraint as an online entry-deterrence strategy.
\newblock {\em Management Science}, 52(11):1799--1809.

\bibitem{zhang2010crafting}
 Zhang, J.,  Farris, P.W.,  Irvin, J.W., Kushwaha, T., Steenburgh, T., 
  and B.~A. Weitz. 2010.
\newblock Crafting integrated multichannel retailing strategies.
\newblock {\em Journal of Interactive Marketing}, 24(2):168--180.

\bibitem{ZhaoZheng2000}
Zhao W., and Y.-S. Zheng. 2000.
\newblock Optimal dynamic pricing for perishable assets with nonhomogeneous
  demand.
\newblock {\em Management Science}, 46(3):375--388.

\end{thebibliography}
\end{document}